\documentclass[12pt]{article}
\input tcilatex
\input xy
\xyoption{all}
\usepackage{amsfonts,epsf}
\newcommand{\psdiag}[3]{\hspace{1mm}\raisebox{-#1mm}{\epsfysize#2mm
\epsffile{#3.eps}}\hspace{1mm}}

\begin{document}
\author{Rui Pedro Carpentier\footnote{rcarpent@math.ist.utl.pt}\\
\\
{\small\it Departamento de Matem\'{a}tica} \\ {\small\it Centro de
An\'{a}lise Matem\'{a}tica, Geometria e Sistemas Din\^{a}micos}\\
{\small\it Instituto Superior T\'{e}cnico}\\ {\small\it Avenida
Rovisco Pais, 1049-001 Lisboa}\\ {\small\it Portugal}}
\title{Three-Colorings of Cubic Graphs and Tensor Operators}
\maketitle

\newpage

\begin{abstract}

Penrose's work \cite{8} established a connection between the edge $3$-colorings of cubic planar graphs and tensor algebras. We exploit this point of view in order to get algebraic representations of the category of cubic graphs with free ends.

\end{abstract}

{\bf keywords:} $3$-colorings of cubic planar graphs, tensor algebras, Penrose invariant, monoidal categories.

\newpage

\section{Introduction}

Although it first appeared as a simple geometric curiosity the
Four Color Problem became one of the most important fields of
research in discrete mathematics linking several areas and having
dozens of equivalent formulations.

The original statement said the following:

\begin{theorem} ($4$-CT) Every planar map can be colored using no more than
four colors in such a way that no pair of adjacent regions receive
the same color.
\end{theorem}

Avoiding the rigorous definitions of map, region or adjacent
regions, this result can be given an equivalent but simpler
statement:

\begin{theorem} Every planar simple graph can be colored using no more than
four colors in such a way that no pair of adjacent vertices
receive the same color.
\end{theorem}

We say that a graph with such a coloring is a $4$-colorable graph
or it has a (vertex) $4$-coloring.

Since any simple planar graph can be embedded in the graph (i.e.
the $1$-skeleton) of a triangulation of the sphere, the $4$-Color
Theorem is equivalent to the following theorem:

\begin{theorem} The graph of any sphere triangulation is
$4$-colorable or has a loop edge.
\end{theorem}

Now regarding a $4$-coloring $\phi$ on the graph of a sphere
triangulation ${\cal T}$ as a $0$-cochain in the simplicial
cohomology of ${\cal T}$ with coefficients in the field of order
$4$, $\mathbb{F}_4$, its coboundary $\delta\phi$ gives a
$3$-coloring on the edges (with colors in $\mathbb{F}_4\setminus
\{0\}$) such that, for any face $f$ of ${\cal T}$, the three edges
$e_1$, $e_2$ and $e_3$ of its boundary receive different colors
(the only way to have
$\delta\phi(e_1)+\delta\phi(e_2)+\delta\phi(e_3)
=\delta\phi(\partial f)=0$ with
$\delta\phi(e_i)\in\mathbb{F}_4\setminus \{0\}$). Note that, since
$\mathbb{F}_4$ is a field of characteristic $2$, it does not
matter what order the simplices (faces, edges or vertices) of
${\cal T}$ have.

On the other hand if we have a $3$-coloring $\psi$ on the edges of
the triangulation ${\cal T}$ assigning different colors to the
three edges of the boundary of any face of ${\cal T}$, then $\psi$
can be regarded as a closed $1$-cochain.
Thus $\psi$ should be the coboundary of some $0$-cochain $\phi$
which would be a (vertex) $4$-coloring of the triangulation ${\cal
T}$.

This proves a result due to Tait that says that the Four Color
Theorem is equivalent to the following proposition:

\begin{theorem} Every planar bridgeless cubic graph is edge $3$-colorable.
\end{theorem}

A {\it cubic graph} is a graph where each vertex is adjacent to
three edges. If a cubic graph is planar then it is the dual graph
of a triangulation of a sphere. A graph is {\it bridgeless} if
there is no edge that after being removed increases the number of
the connected components. A planar cubic graph is bridgeless if
and only if it is the dual graph of a triangulation without loops
of a sphere.

Much of the research in this area focuses mainly on the Tait
version of the Four Color Theorem. The references \cite{6,7} provied a good overview about the Four Color Theorem and its ramifications.

\section{Category of cubic graphs (with free ends)}

It is possible to study the edge $3$-colorings of cubic graphs by
introducing a category of cubic graphs with free ends.

Consider the following (monoidal) category $\mathbf{CG}$. The
objects of $\mathbf{CG}$ are the non-negative integer numbers and
a morphism from $m$ to $n$ is a regular immersion of a cubic graph
with $m+n$ free ends  in the strip $\mathbb{R}\times [0,1]$ such
that the free ends are placed at the points $(1,1)$, ... ,$(m,1)$
and $(1,0)$, ... ,$(n,0)$ (see the next figure).

$$\psdiag{10}{20}{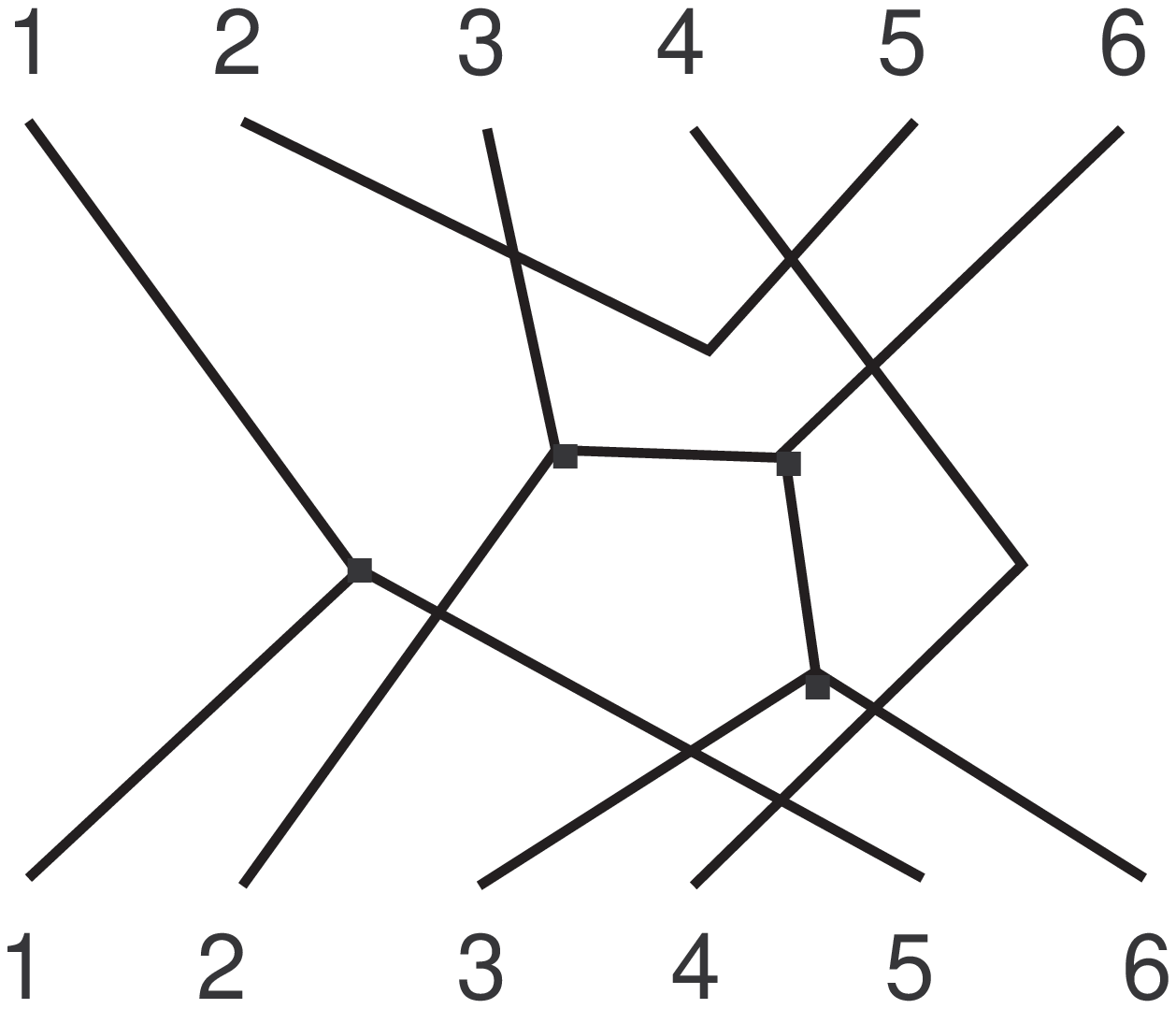}$$

To simplify the treatment we consider piecewise linear immersions
rather than smooth immersions.

The composition in this category is defined in the following way.
Given two morphisms $g_1 :l\rightarrow m$ and $g_2 :m\rightarrow
n$  then their composition $g_2 g_1 :l\rightarrow n$ would be the
immersion $g(f(g_1)\cup g_2)$ where $f(x,y)=(x,y+1)$ and
$g(x,y)=(x,y/2)$ (see the next figure)\footnote{In this paper, the downward direction composition is used, some
authors use the opposite direction.}.

$$g_1 = \psdiag{8}{16}{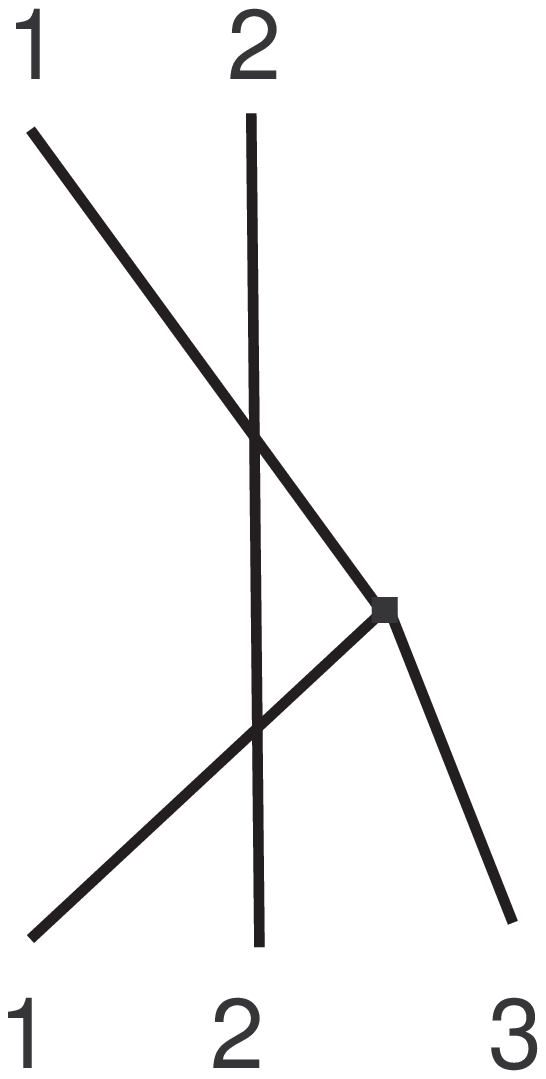}\quad \mbox{and} \quad g_2 =
\psdiag{8}{16}{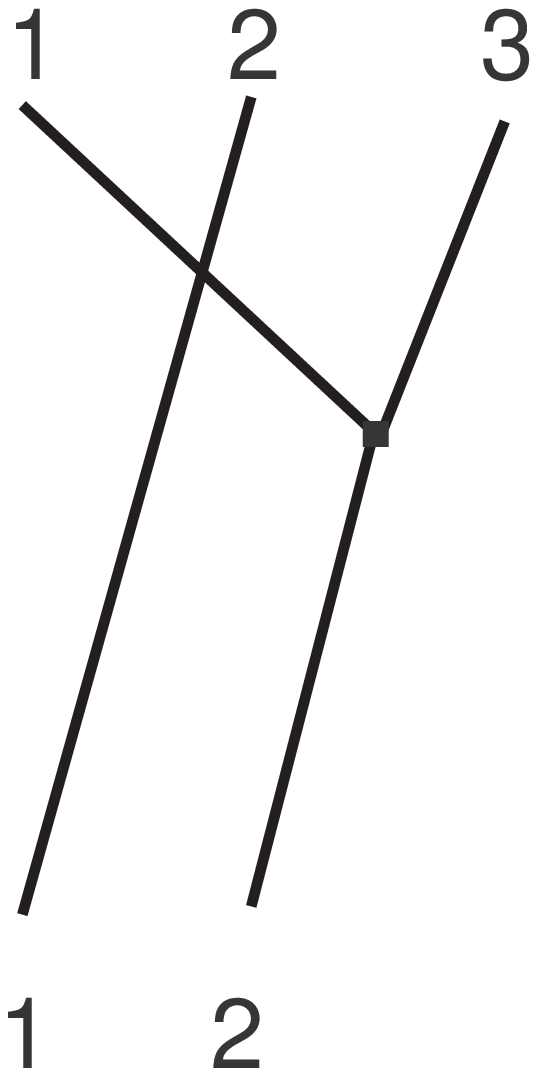}
 \longrightarrow g_2 g_1 =
\psdiag{8}{16}{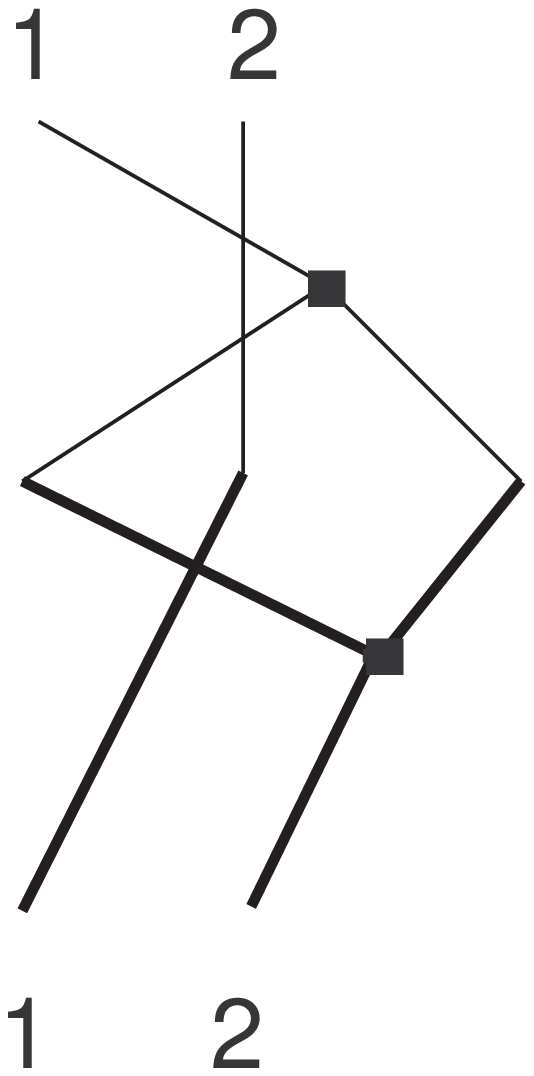}$$

This category places the same role for cubic graphs as the
category of the tangles is for links and like the latter it has a
monoidal structure. Given two morphisms $g_1 :k\rightarrow l$ and
$g_2 :m\rightarrow n$  we get a new morphism $g_1\otimes g_2:
k+m\rightarrow l+n$ by putting the two graph immersions side by
side (see the next figure).

$$g_1 = \psdiag{8}{16}{g1}\quad \mbox{and} \quad g_2 =
\psdiag{8}{16}{g2}
 \longrightarrow g_2\otimes g_1 =
\psdiag{8}{16}{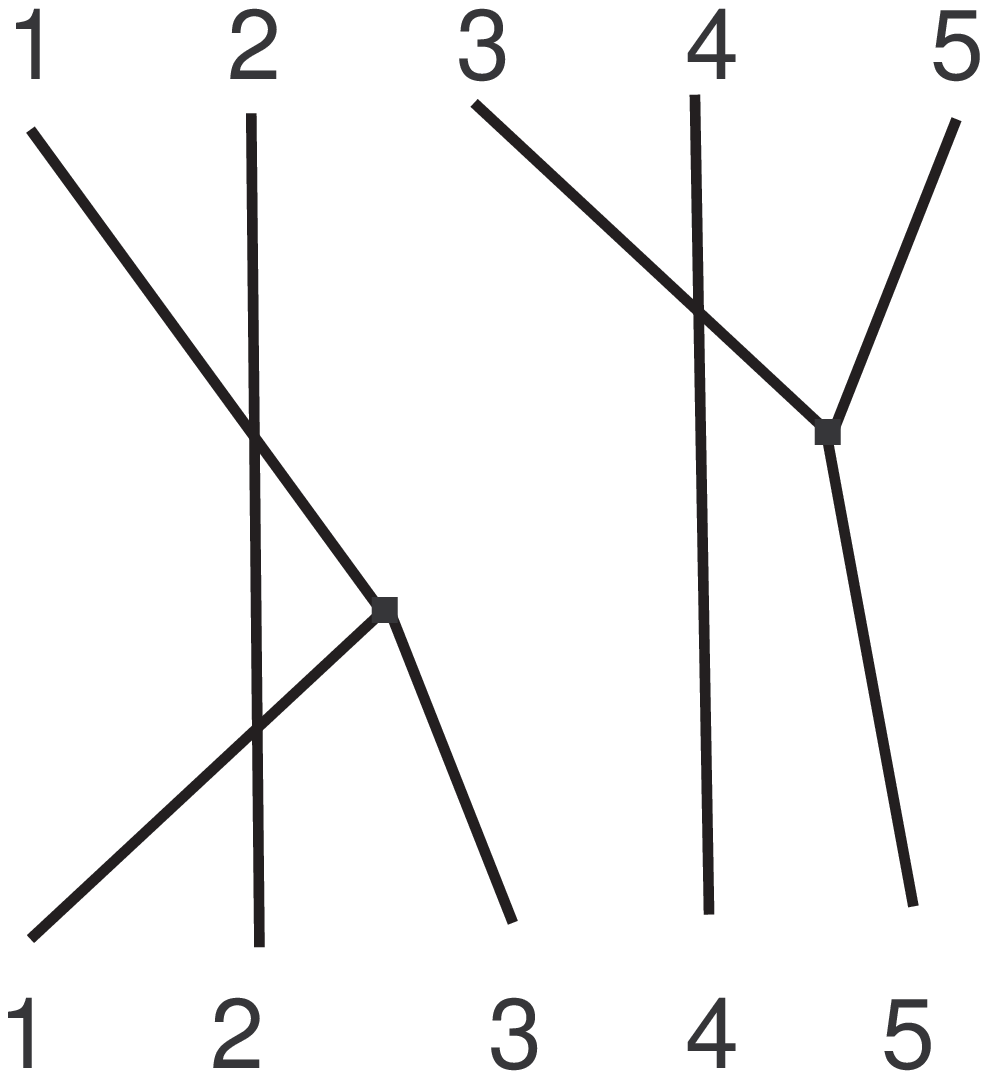}$$

It is easy to see that with these two operations the category
${\cal CG}$ is generated by the following morphisms:

$$\psdiag{4}{8}{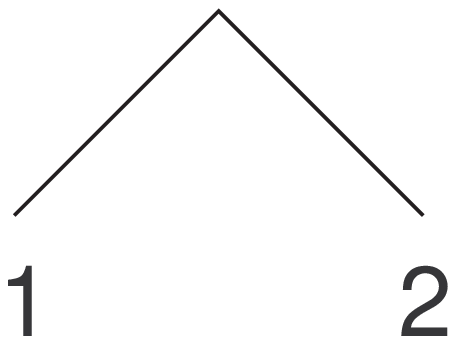}\quad\cap:0\rightarrow 2\quad  ,\quad \cup
:2\rightarrow 0 \quad\psdiag{4}{8}{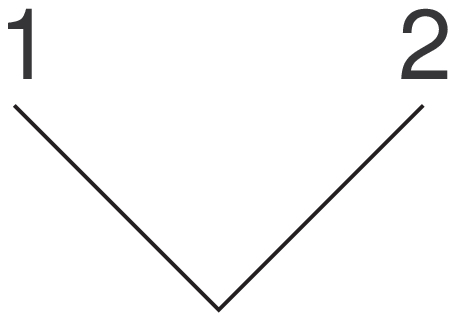}$$

$$\psdiag{4}{8}{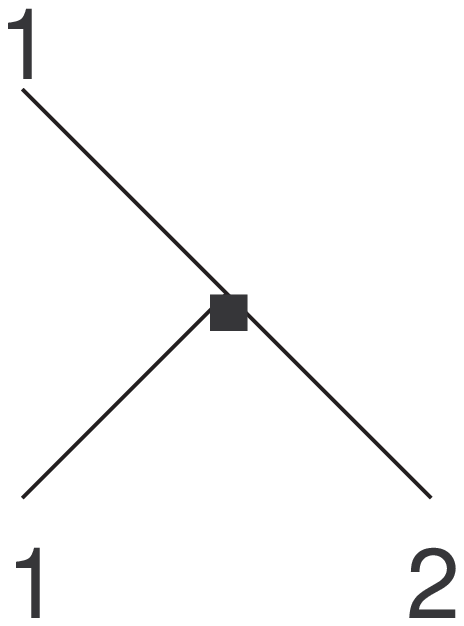}\quad\lambda:1\rightarrow 2\quad \quad ,\quad y
:2\rightarrow 1 \quad\psdiag{4}{8}{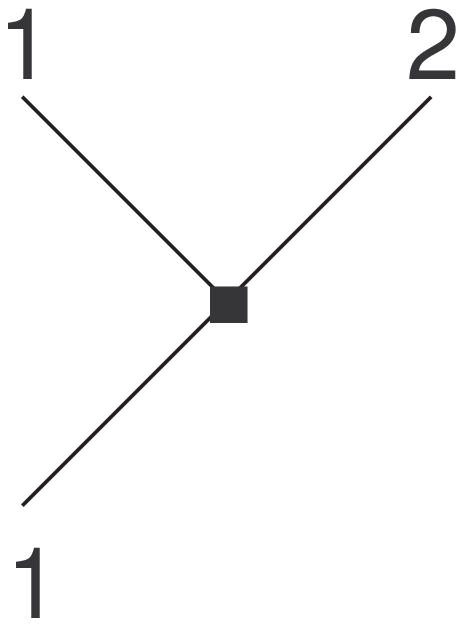}$$

$$\psdiag{4}{8}{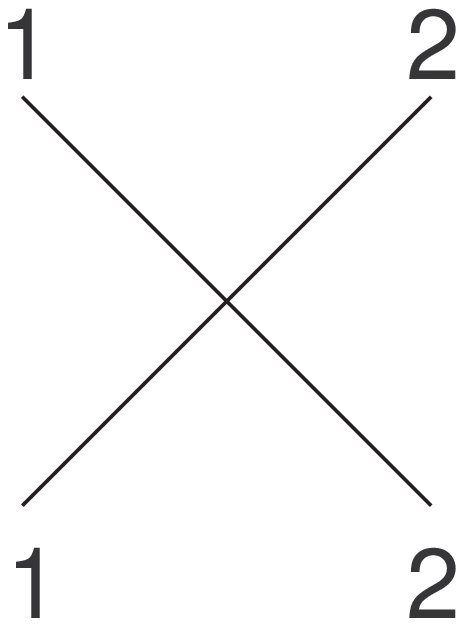}\quad x:2\rightarrow 2\quad \quad \mbox{and}
\quad I:=id_1 :1\rightarrow 1 \quad\psdiag{4}{8}{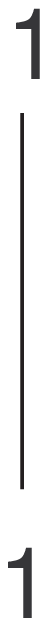}$$

There are some relations that these generators should satisfy:

$$(\cup\otimes I)(I\otimes \cap)=I=(I\otimes \cup)(\cap\otimes
I)\quad\psdiag{4}{8}{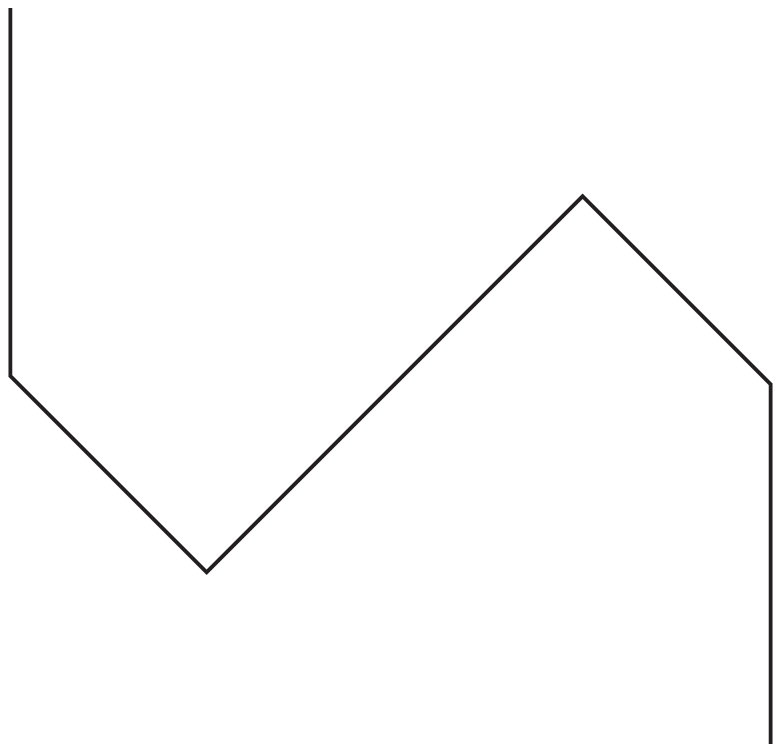}=\psdiag{4}{8}{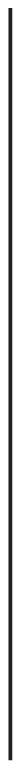}=\psdiag{4}{8}{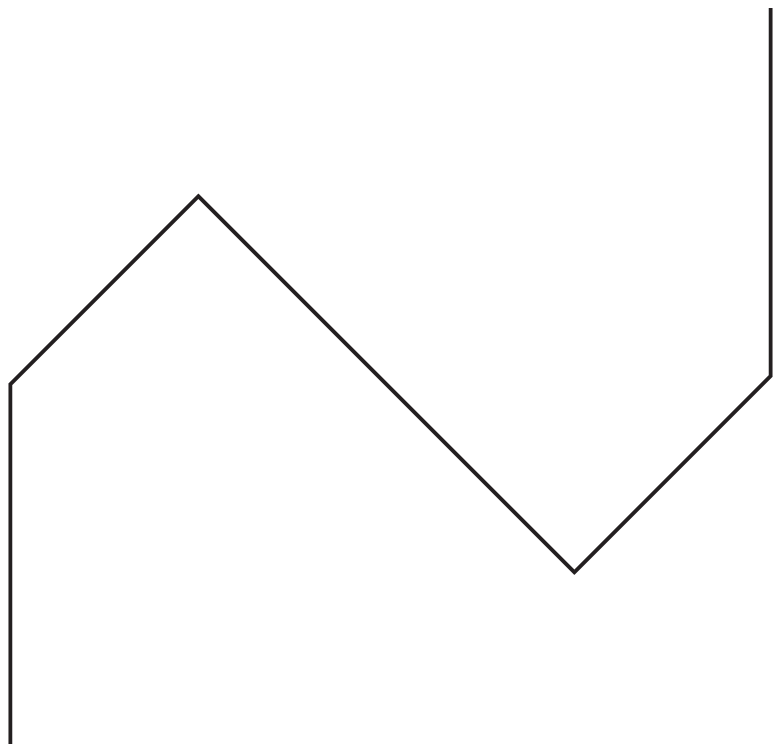}$$

$$(x\otimes I)(I\otimes \cap)=(I\otimes x)(\cap\otimes
I)\quad\psdiag{4}{8}{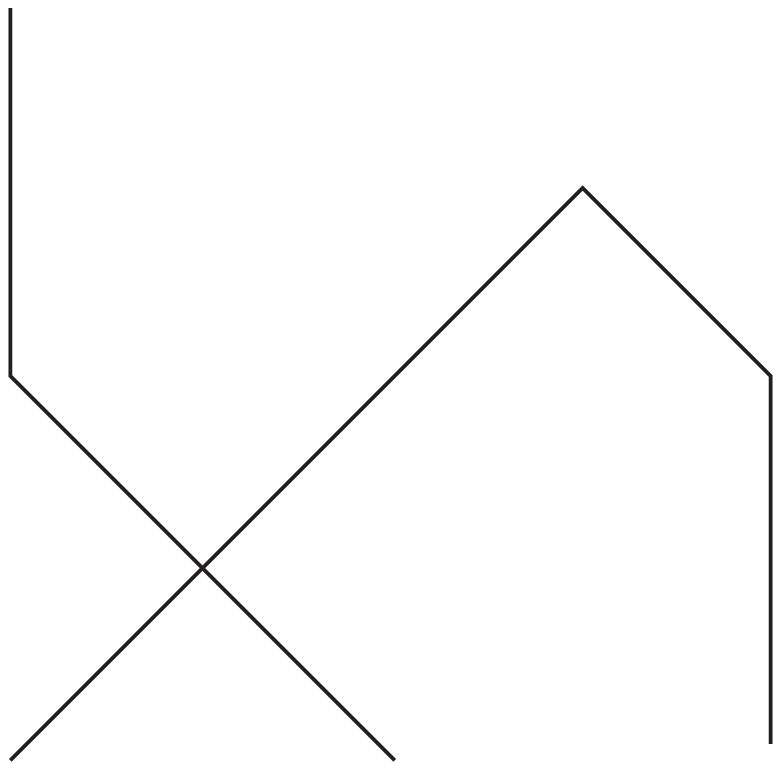}=\psdiag{4}{8}{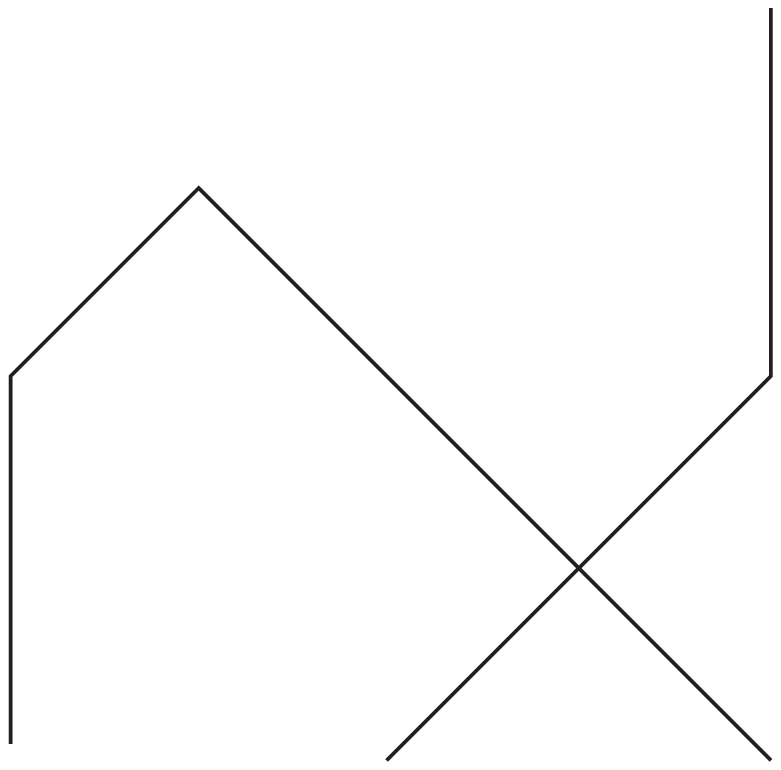}$$

$$\cup x=\cup\quad\psdiag{4}{8}{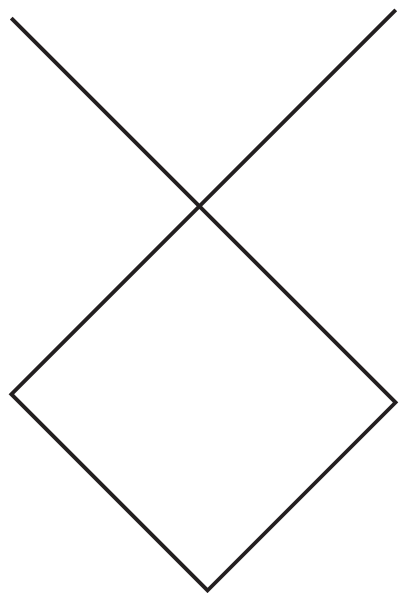}=\psdiag{4}{8}{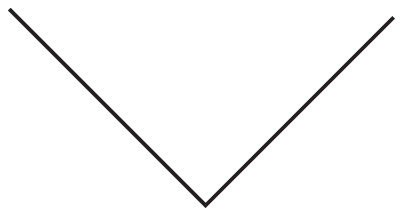}$$

$$xx=I\otimes I\quad\psdiag{4}{8}{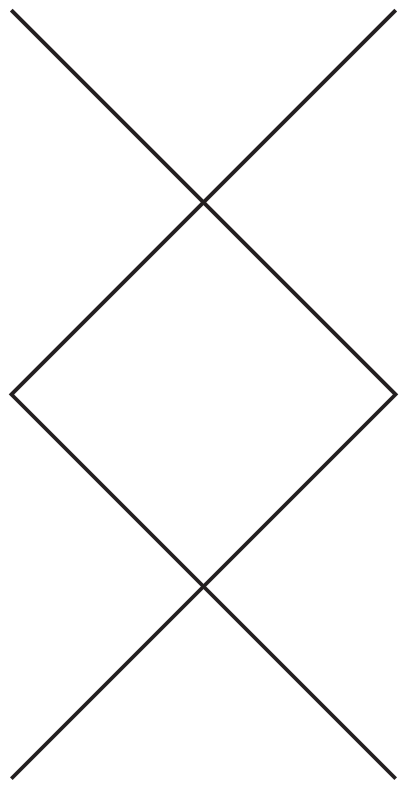}=\psdiag{4}{8}{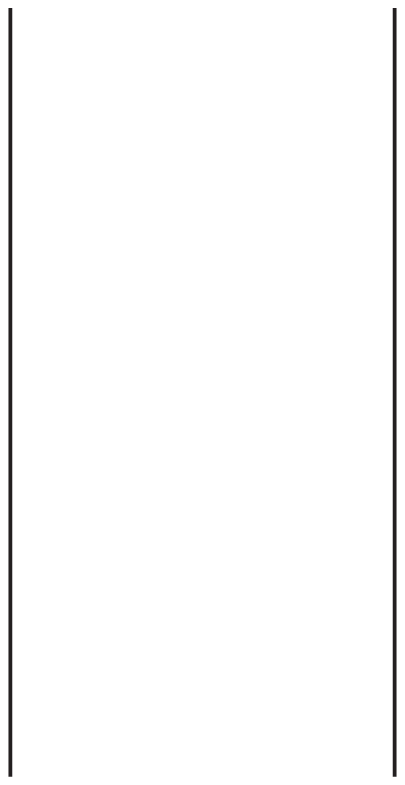}$$

$$(I\otimes x)(x\otimes I)(I\otimes x)=(x\otimes I)(I\otimes
x)(x\otimes I)\quad\psdiag{4}{8}{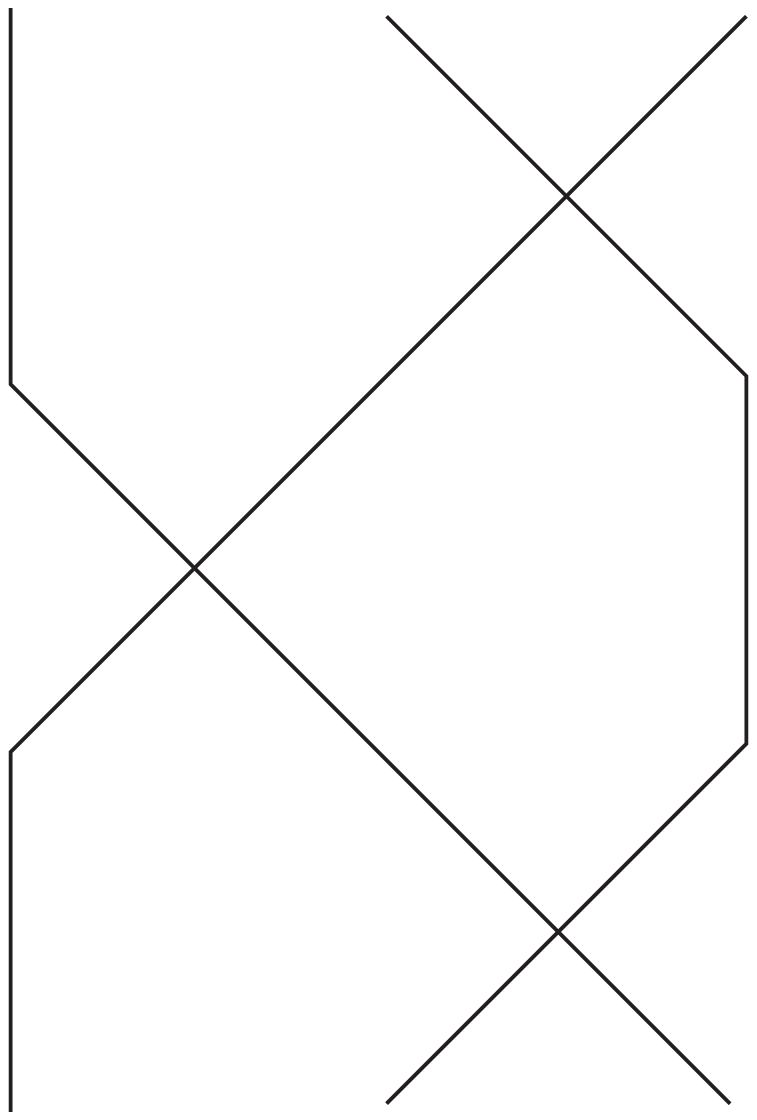}=\psdiag{4}{8}{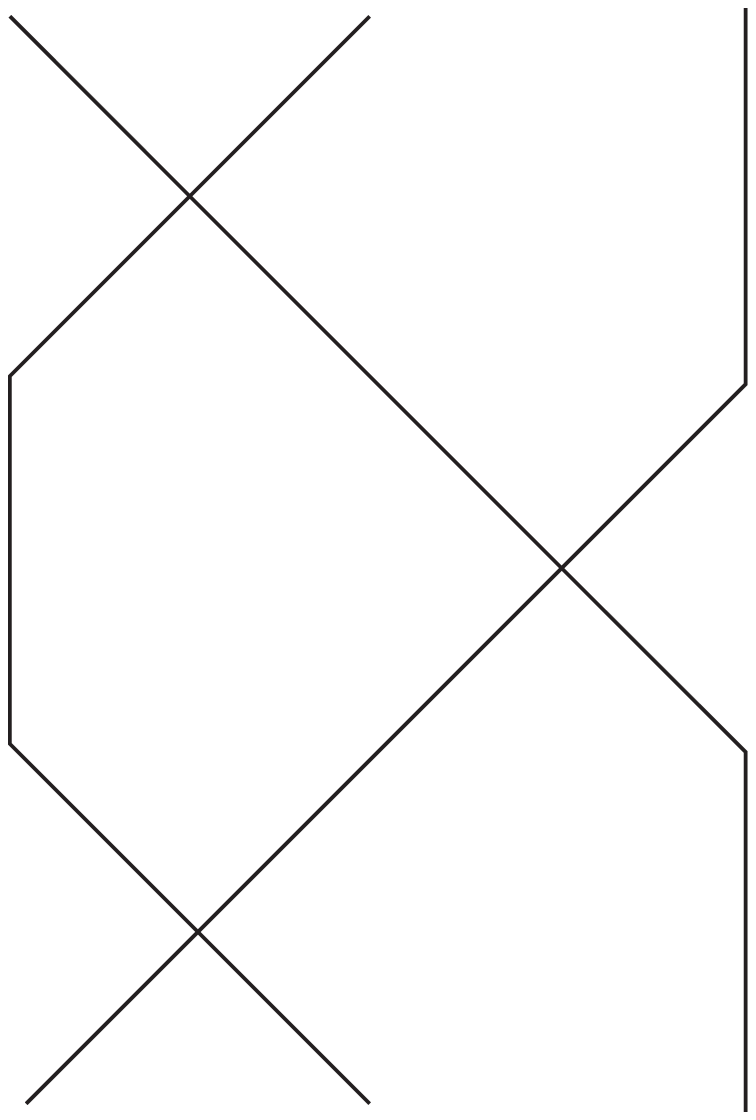}$$

$$(\lambda\otimes I)\cap=I=(I\otimes
\lambda)\cap\quad\psdiag{4}{8}{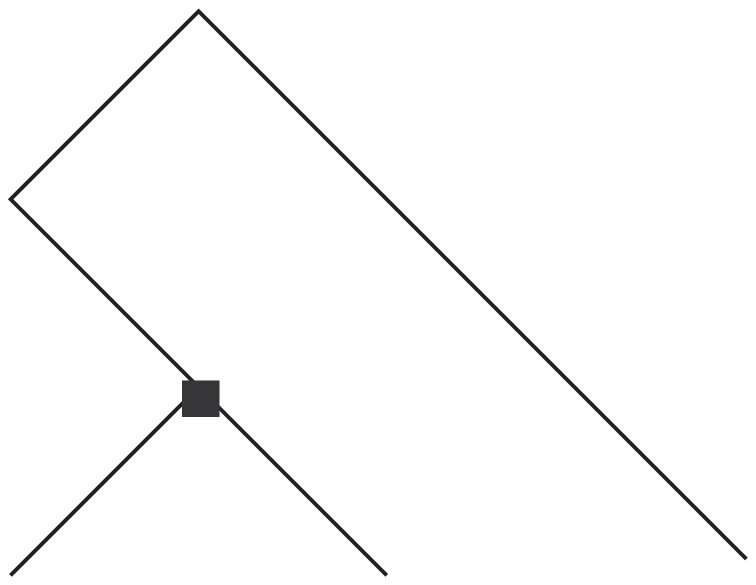}=\psdiag{4}{8}{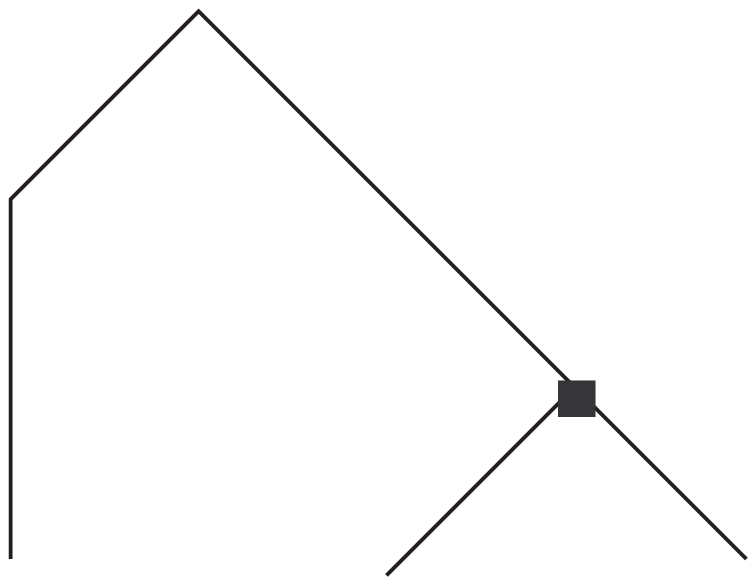}$$

$$(y\otimes I)(I\otimes \cap)=\lambda=(I\otimes y)(\cap\otimes
I)\quad\psdiag{4}{8}{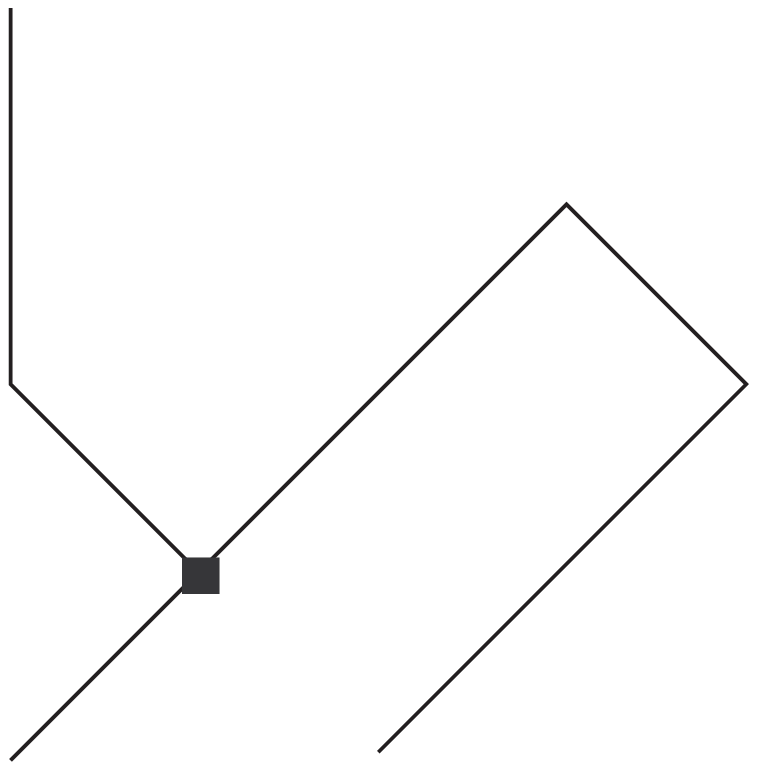}=\psdiag{4}{8}{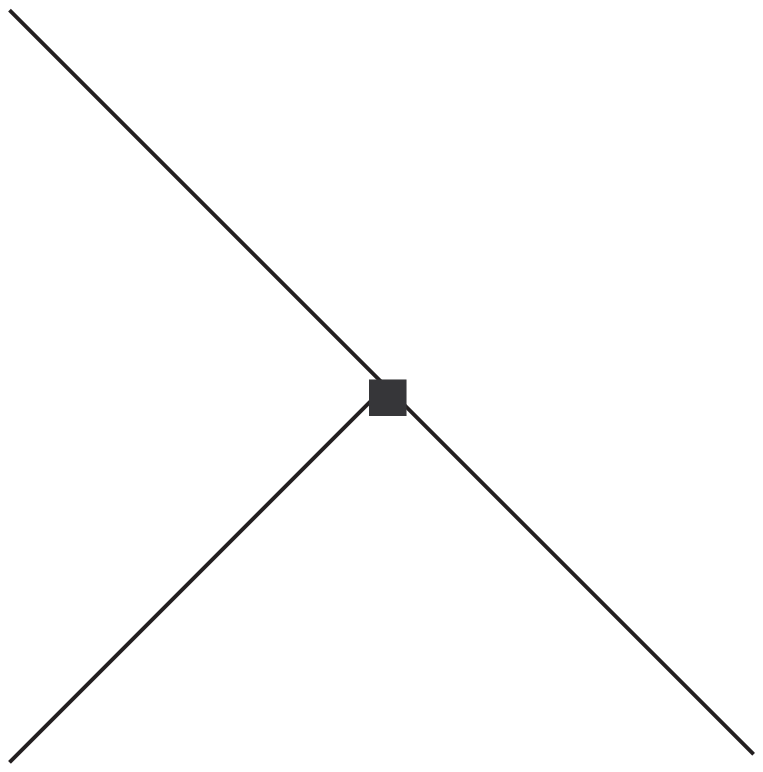}=\psdiag{4}{8}{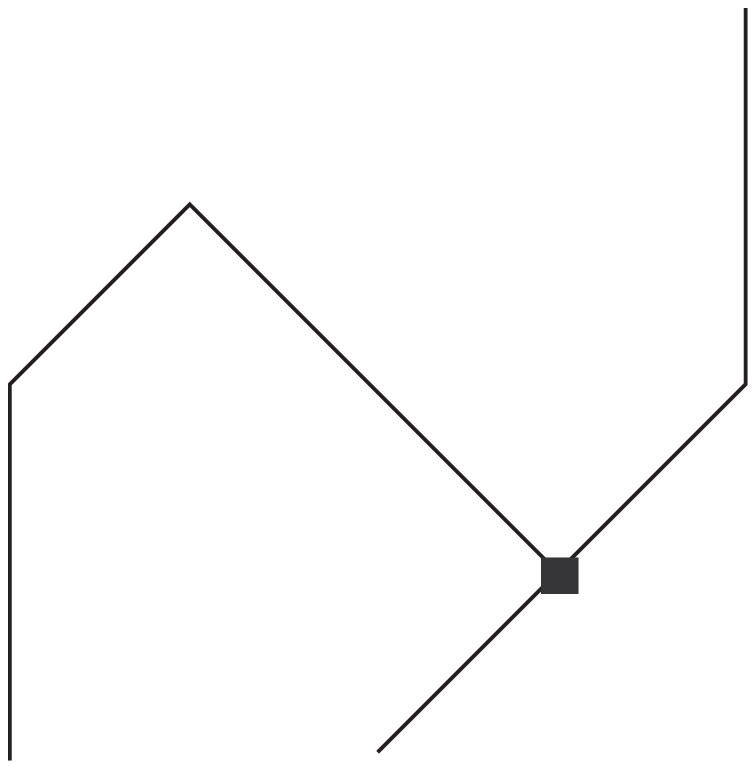}$$

$$(x\otimes I)(I\otimes \lambda) x=(I\otimes x)(\lambda\otimes
I)\quad\psdiag{4}{8}{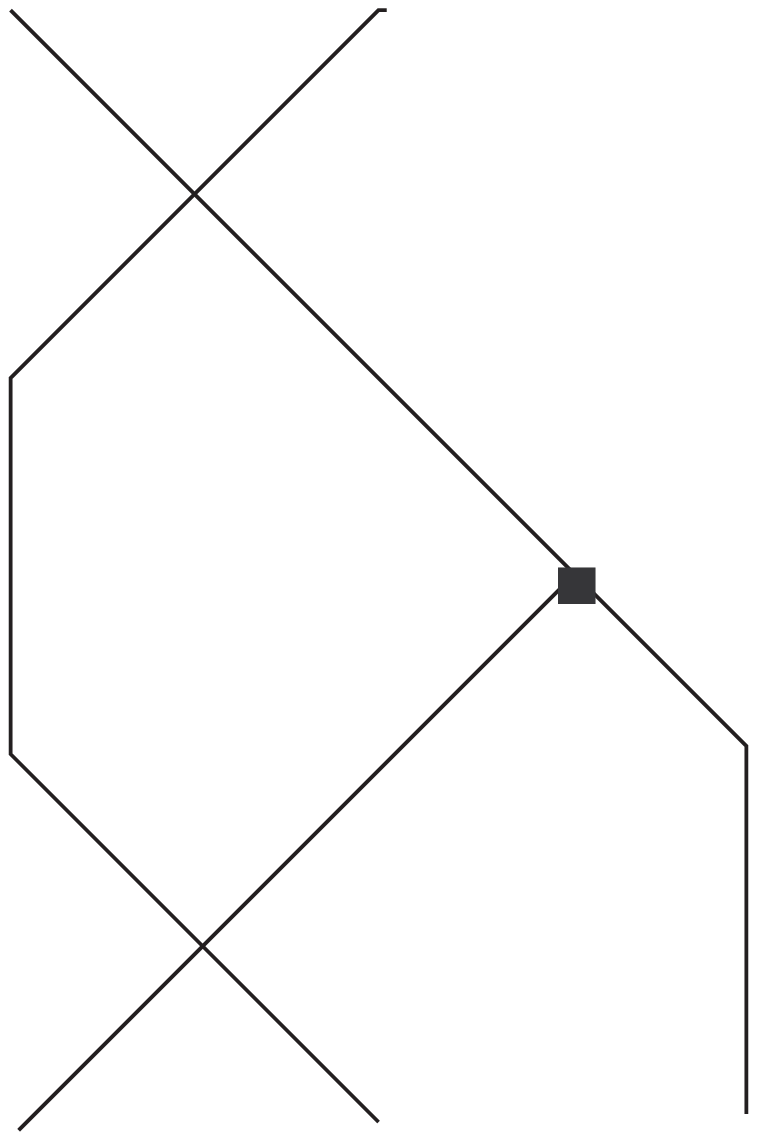}=\psdiag{4}{8}{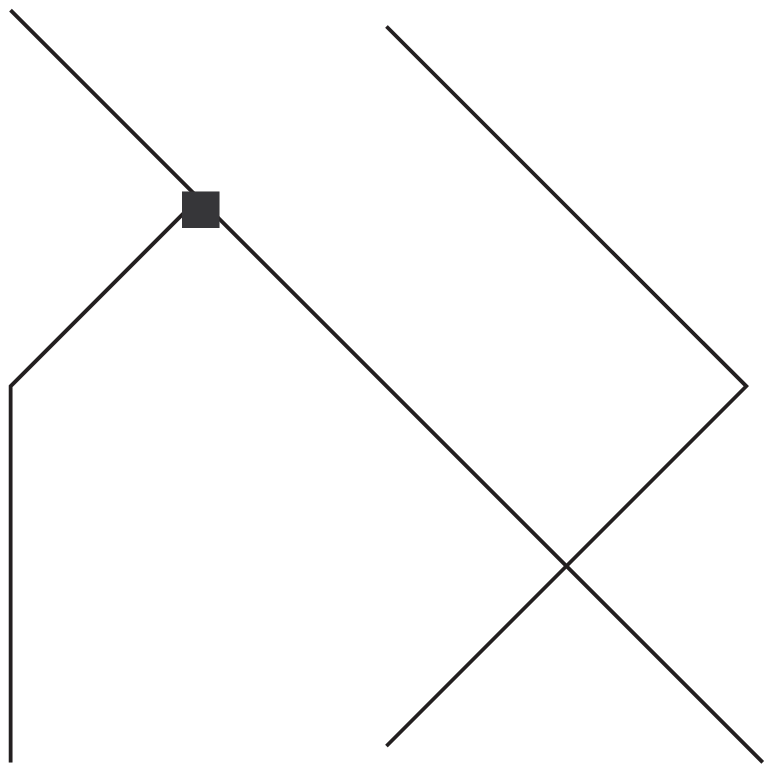}$$

Besides these relations, as a strict monoidal category, $\mathbf{CG}$
should satisfy the following equality:

Given two morphisms $f:k\rightarrow l$ and $g:m\rightarrow n$

$$f\otimes g=(f\otimes id_n)(id_k \otimes g)=(id_l \otimes
g)(f\otimes id_m)$$

As a consequence of this identity and the previous relations we have that this representation is invariant under ambient isotopies.


If we drop the generator $x$ we get a subcategory $\mathbf{PCG}$
of $\mathbf{CG}$ which only contains planar cubic graphs with free
ends.

Now let $\mathbb{K}$ be a field of characteristic zero and let $V$
be a $3$-dimensional $\mathbb{K}$-vector space. We fix a canonical
basis $\{e_1,e_2,e_3\}$ for $V$ and introduce the following
(monoidal) functor from $\mathbf{CG}$ to the (monoidal) category
$\mathbf{Vect}_\mathbb{K}$ of vector spaces over $\mathbb{K}$:

$$F:\mathbf{CG}\rightarrow \mathbf{Vect}_\mathbb{K}$$ defined on
the objects by $$F(n)=V^{\otimes n} \quad (F(0)=\mathbb{K})$$ and
on the morphisms by the following definitions on the generators:

$$\begin{array}{cccc}
  F(\cap): & \mathbb{K} & \longrightarrow & V\otimes V \\
   & \alpha & \longmapsto & \alpha e_i\otimes e_i
\end{array}$$

$$\begin{array}{cccc}
  F(\cup): & V\otimes V & \longrightarrow & \mathbb{K} \\
   & \sum_{i,j=1}^{3}\alpha_{i,j} e_i\otimes e_j & \longmapsto &
   \sum_{i=1}^{3}\alpha_{i,i}
\end{array}$$

$$\begin{array}{cccc}
  F(\lambda): & V & \longrightarrow & V\otimes V \\
   & \sum_{i=1}^{3}\alpha_{i}e_i & \longmapsto &
   \sum_{\{i,j,k\}=\{1,2,3\}}\alpha_{i} e_j\otimes e_k
\end{array}$$

$$\begin{array}{cccc}
  F(y): & V\otimes V & \longrightarrow & V \\
   & \sum_{i,j=1}^{3}\alpha_{i,j}e_i\otimes e_j & \longmapsto &
   \sum_{\{i,j,k\}=\{1,2,3\}}\alpha_{i,j} e_k
\end{array}$$

$$\begin{array}{cccc}
  F(x): & V\otimes V & \longrightarrow & V\otimes V \\
   & \sum_{i,j=1}^{3}\alpha_{i,j}e_i\otimes e_j & \longmapsto &
   \sum_{i,j=1}^{3}\alpha_{i,j}e_j\otimes e_i
\end{array}$$

It is easy to see that this functor is well defined under the
relation (i.e. $F[(\cup\otimes I)(I\otimes \cap)]=F(I)=F[(I\otimes
\cup)(\cap\otimes I)]$, $F[(x\otimes I)(I\otimes
\cap)]=F[(I\otimes x)(\cap\otimes I)]$, ..., etc). In fact, we
have the following theorem.

\begin{theorem} If a morphism $g:m\rightarrow n$ represents a
cubic graph with $m+n$ free ends then for each element,
$e_{i(1)}\otimes\cdots\otimes e_{i(m)}$, of the canonical basis of
$V^{\otimes m}$ $$F(g)(e_{i(1)}\otimes\cdots\otimes e_{i(m)})=
\sum_{j(1),...,j(n)} \chi_{j(1),...,j(n)}^{i(1),...,i(m)}
e_{j(1)}\otimes\cdots\otimes e_{j(n)}$$ where
$\chi_{j(1),...,j(n)}^{i(1),...,i(m)}$ is the number of edge
$3$-colorings of the graph such that it has the free edges on the
top colored by $i(1),...,i(m)$ (in this order) and the free edges
on the bottom colored by $j(1),...,j(n)$ in this order.
\end{theorem}

\begin{example}

$g:2\rightarrow 2$ : $\psdiag{4}{8}{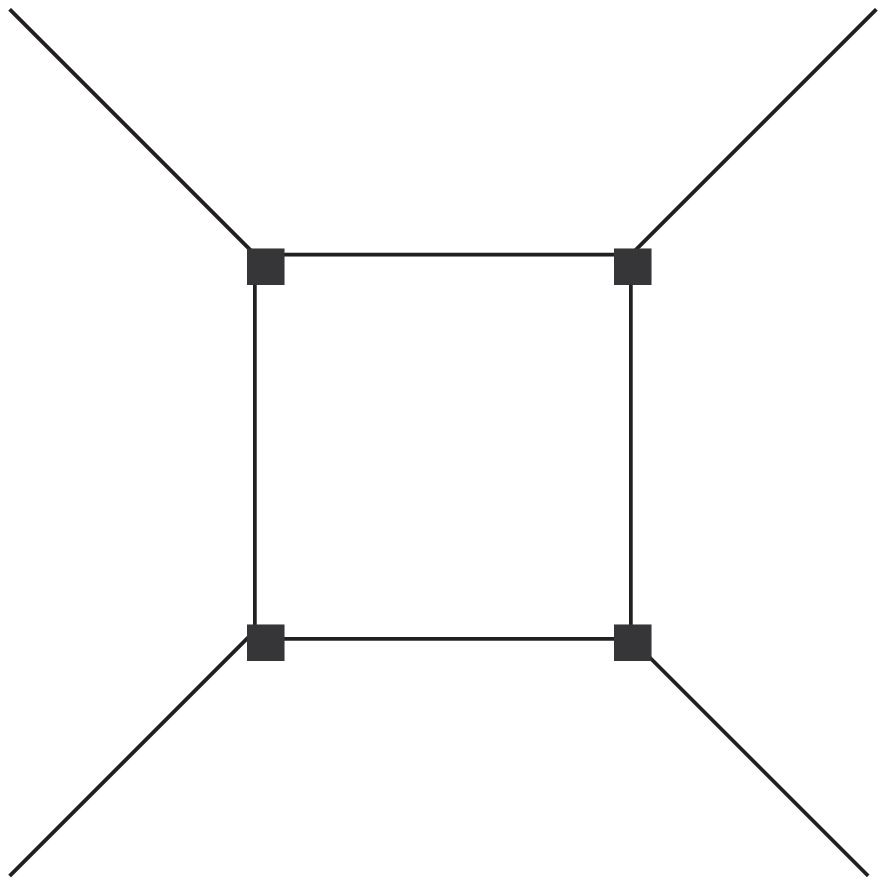}$

egde $3$-colorings with colors $1,1$ on the top:
\{$\psdiag{4}{8}{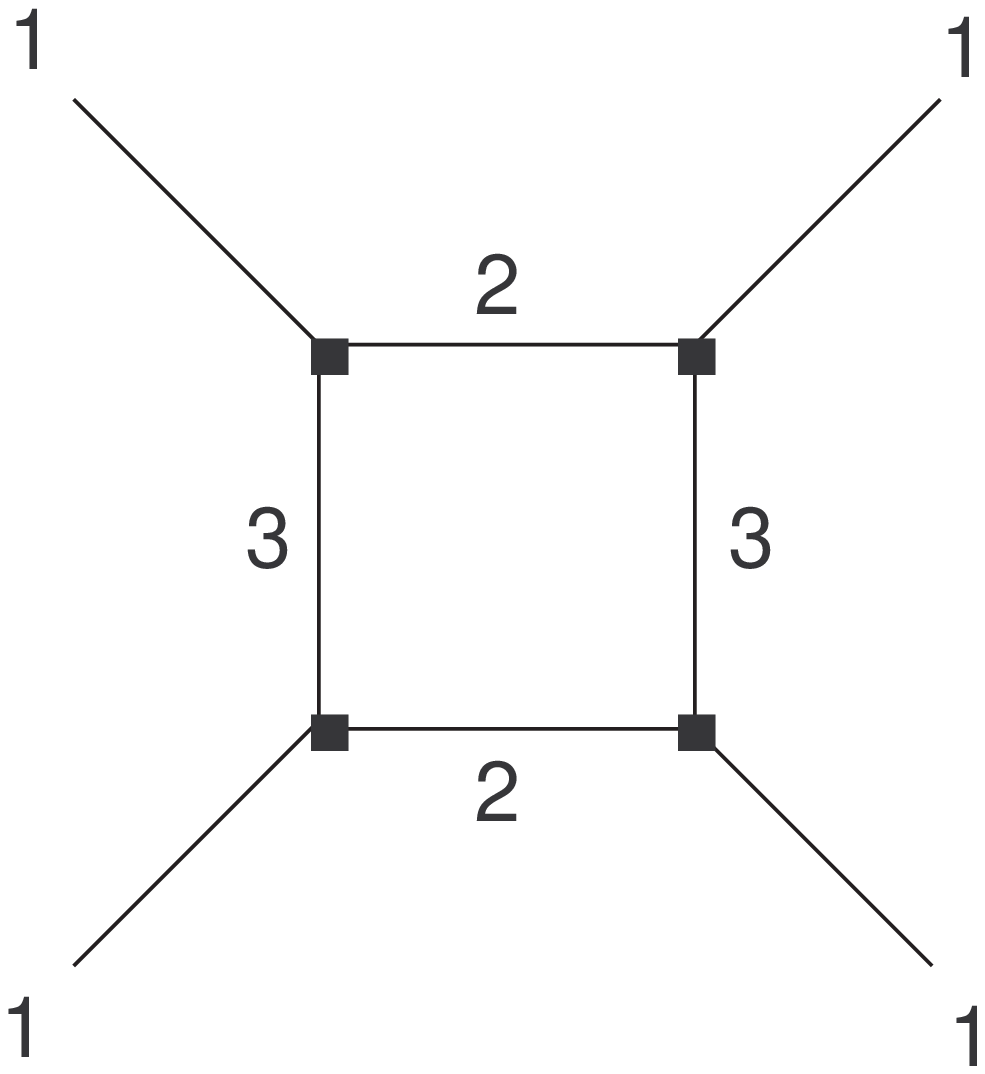}$, $\psdiag{4}{8}{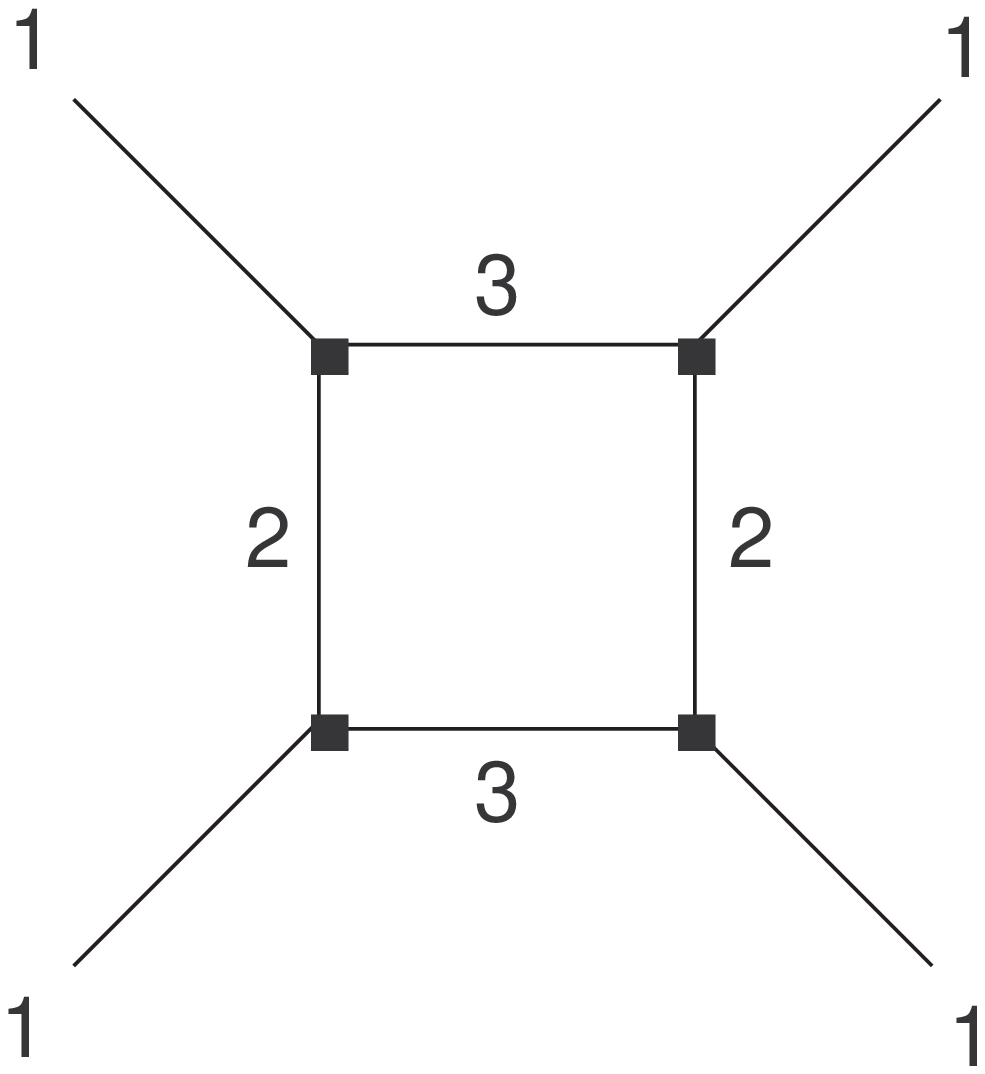}$\},
$\psdiag{4}{8}{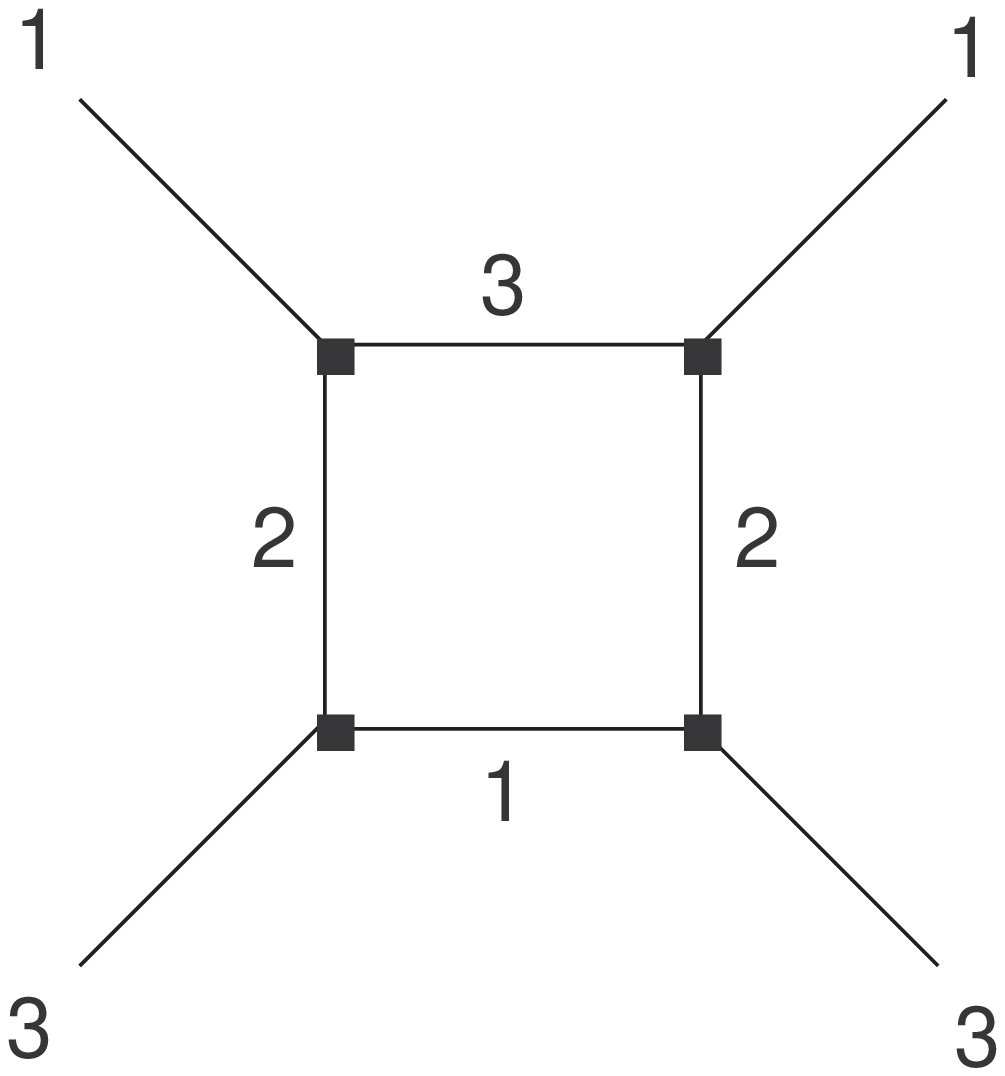}$ and $\psdiag{4}{8}{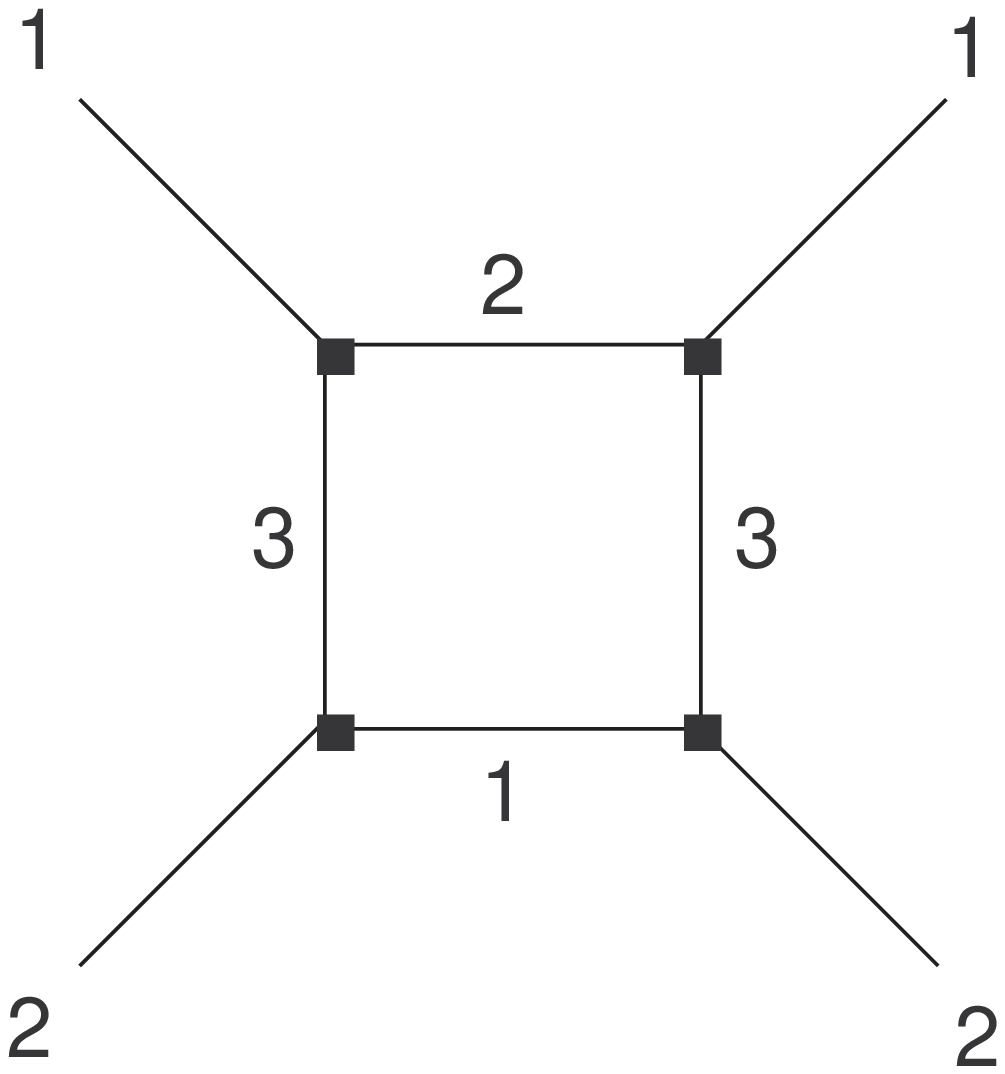}$

$F(g)(e_1\otimes e_1)=2e_1\otimes e_1+e_2\otimes e_2+e_3\otimes
e_3$

\end{example}

\TeXButton{Proof}{\proof}

It is enough to check the statement on the generators $\cap$,
$\cup$, $\lambda$, $y$ and $x$ (which is straightforward) and to
note that the composition and the monoidal operation on
$\mathbf{CG}$ satisfy also the required:

$$\chi_{j(1),...,j(n)}^{i(1),...,i(k)}(gf)=
\sum_{\alpha(1),...,\alpha(m)}\chi_{j(1),...,j(n)}^{\alpha(1),...,\alpha(m)}(g)
\chi_{\alpha(1),...,\alpha(m)}^{i(1),...,i(k)}(f)$$

$$\chi_{j(1),...,j(l+n)}^{i(1),...,i(k+m)}(g\otimes f)=
\chi_{j(1),...,j(l)}^{i(1),...,i(k)}(g)
\chi_{j(l+1),...,j(l+n)}^{i(k+1),...,i(k+m)}(f) $$

\TeXButton{End Proof}{\endproof}

The following corollary is an immediate consequence of the
theorem.

\begin{corollary}
Given a cubic graph $g$ (without free end edges) viewed as a
morphism $g:0\longrightarrow 0$ the value $F(g)(1)$ is the number
of edge $3$-colorings of the graph. In particular $g$ is edge
$3$-colorable if and only if $F(g)(1)\not=0$.
\end{corollary}

Next, we introduce another functor $\tilde{F}:\mathbf{CG}
\longrightarrow \mathbf{Vect}_\mathbb{K}$ which is a small
modification of the functor $F$.

$\tilde{F}$ is equal to $F$ on the objects and on all the
generator morphisms except on the morphism $x$ where

$$\begin{array}{cccc}
  \tilde{F}(x): & V\otimes V & \longrightarrow & V\otimes V \\
   & \sum_{i,j=1}^{3}\alpha_{i,j}e_i\otimes e_j & \longmapsto &
   \sum_{i,j=1}^{3}\tilde{\delta}_{i,j}\alpha_{i,j}e_j\otimes e_i
\end{array}$$

with $\tilde{\delta}_{i,j}=-1+2\delta_{i,j}$ where $\delta_{i,j}$
is the Kronecker delta.

We have that $\tilde{F}=F$ when restricted to the subcategory
$\mathbf{PCG}$ (the planar cubic graphs).

The special feature of the functor $\tilde{F}$ is that it
satisfies the Penrose formula:

$$\tilde{F}(\psdiag{2}{6}{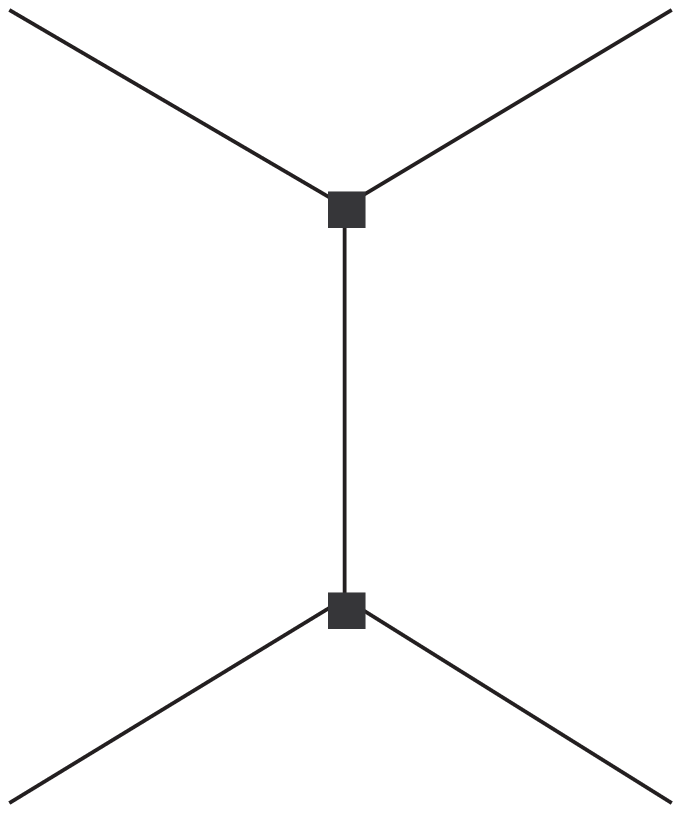})=\tilde{F}(\psdiag{2}{6}{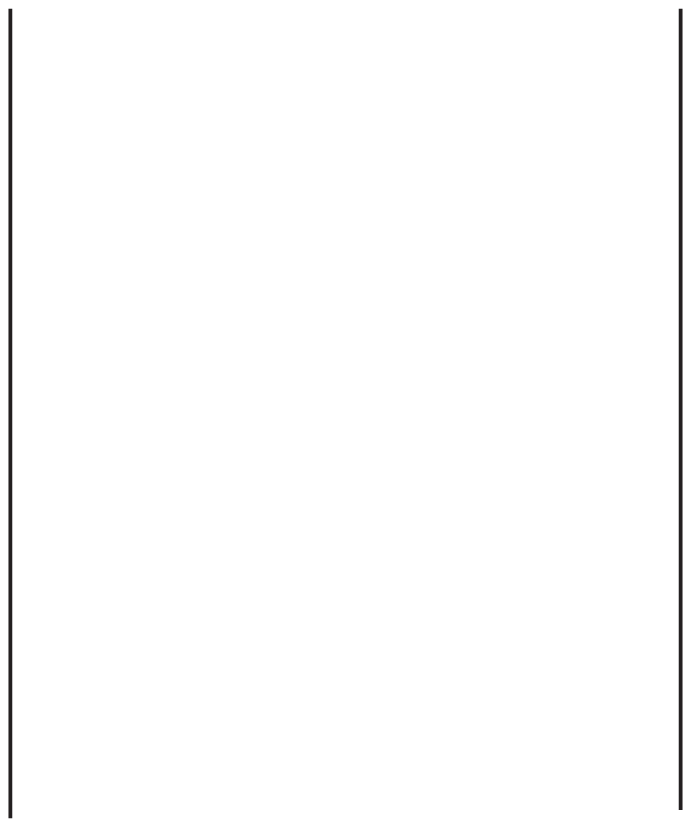})-\tilde{F}(\psdiag{2}{6}{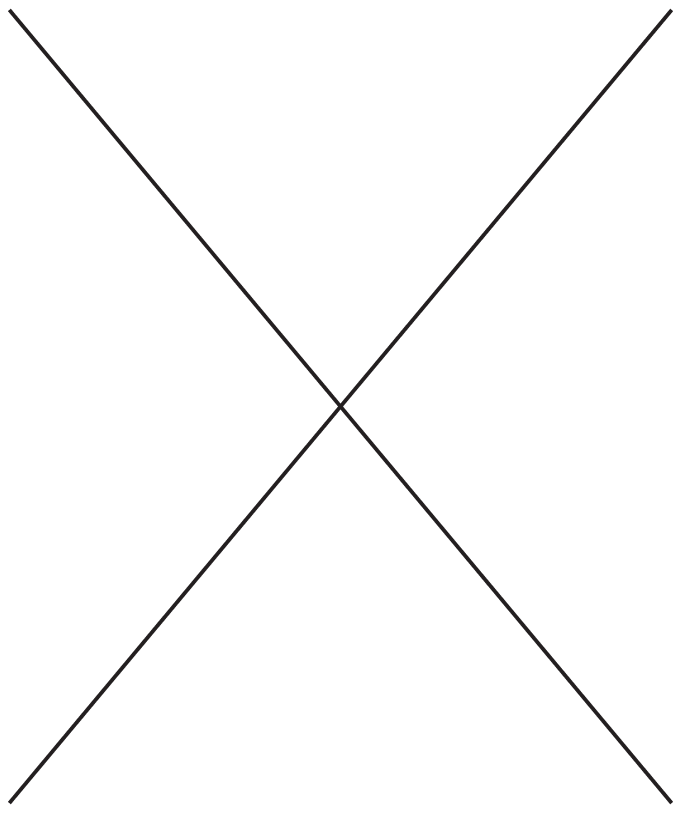})$$

Thus it also satisfies the IHX identity on chinese characters (see
\cite{1})

$$\tilde{F}(\psdiag{2}{6}{pen1})=\tilde{F}(\psdiag{2}{6}{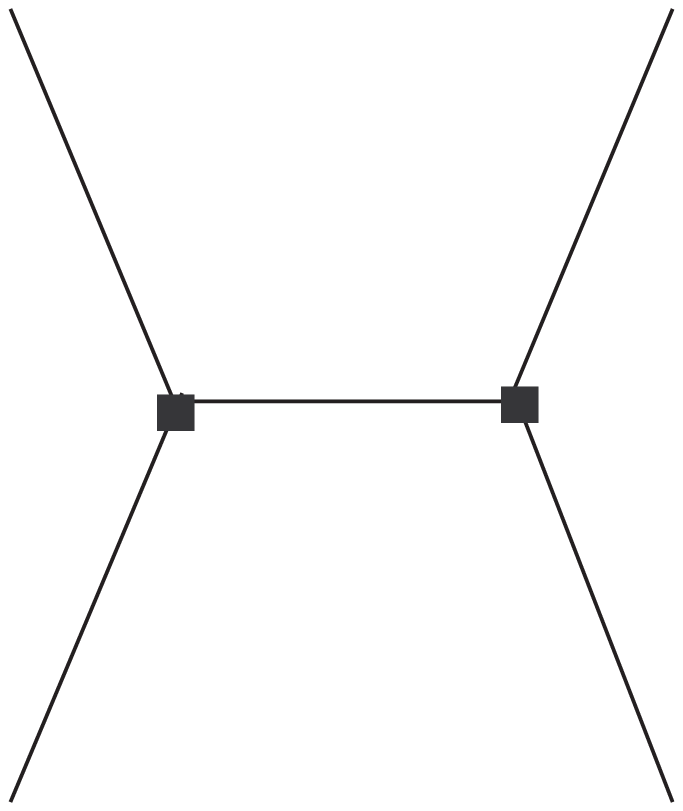})-\tilde{F}(\psdiag{2}{6}{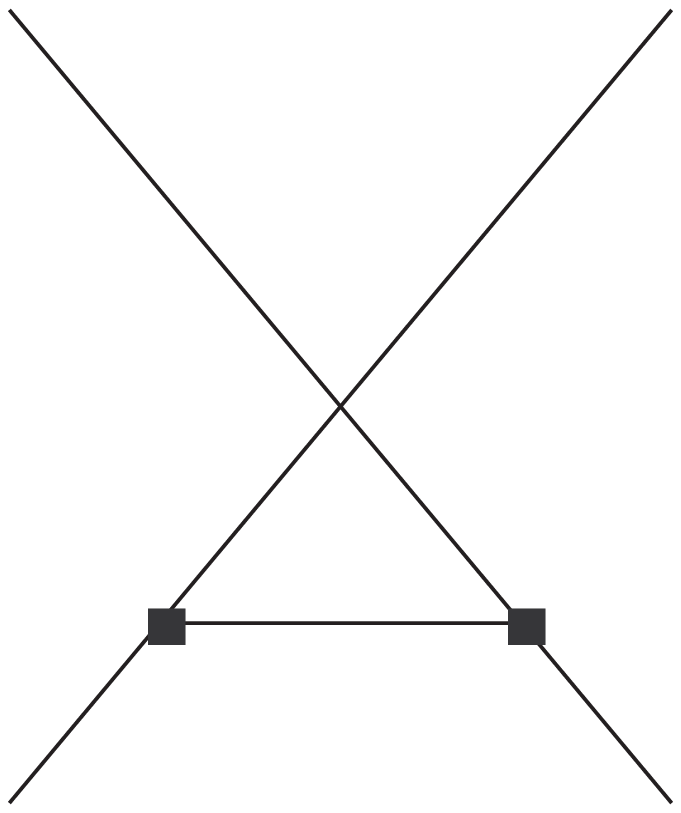})$$

We also have the formula
$$\tilde{F}(\psdiag{2}{6}{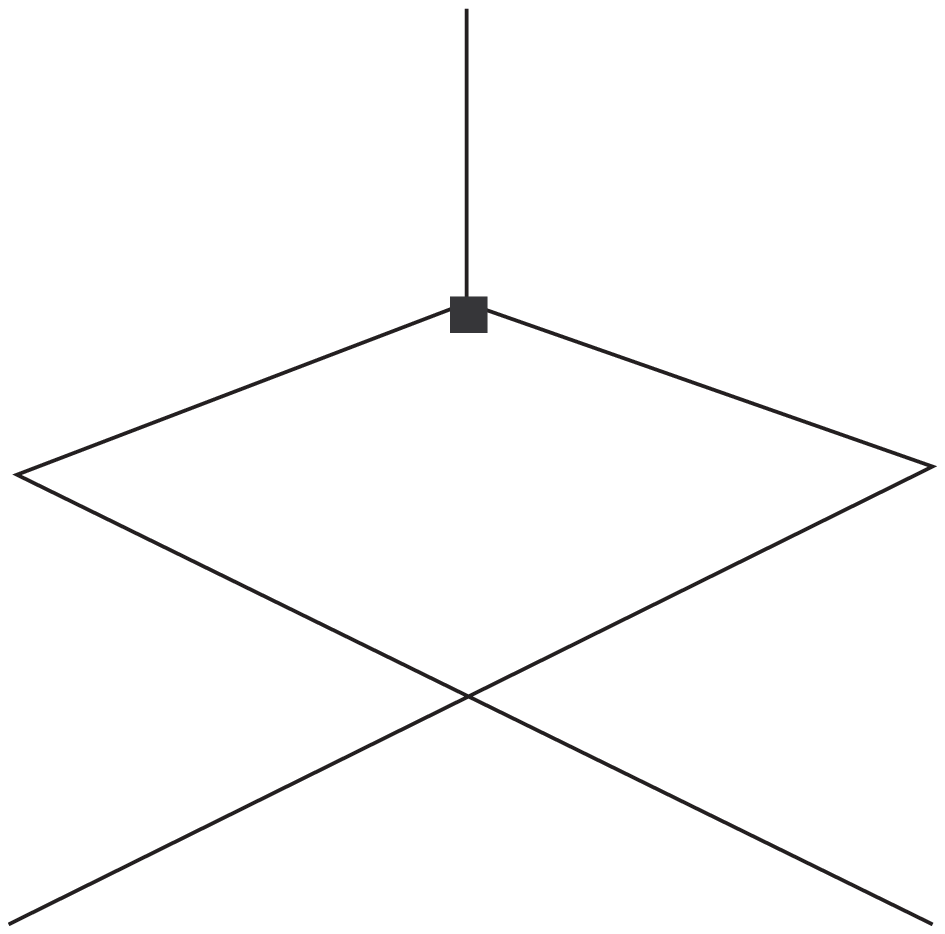})=-\tilde{F}(\psdiag{2}{6}{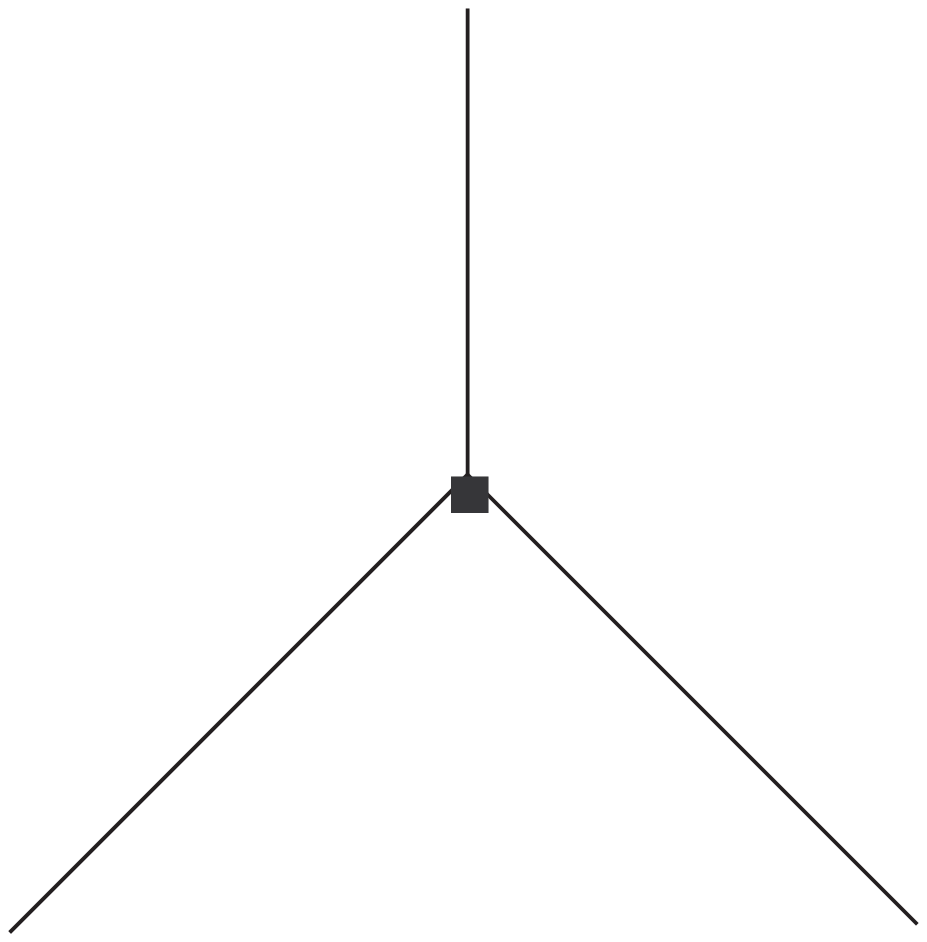})$$

To see what the functor $\tilde{F}$ gives let us introduce the
notion of the {\it sign} of an edge $3$-coloring of a cubic graph.
The sign of an edge $3$-coloring of a cubic graph projected on the
plane is $+1$ or $-1$ if the number of crossings of edges of
different colors is even or odd. Then we have the following
result.

\begin{theorem} If a morphism $g:m\rightarrow n$ represents a
cubic graph with $m+n$ free ends then for each element,
$e_{i(1)}\otimes\cdots\otimes e_{i(m)}$, of the canonical basis of
$V^{\otimes m}$ $$\tilde{F}(g)(e_{i(1)}\otimes\cdots\otimes
e_{i(m)})= \sum_{j(1),...,j(n)}
\tilde{\chi}_{j(1),...,j(n)}^{i(1),...,i(m)}
e_{j(1)}\otimes\cdots\otimes e_{j(n)}$$ where
$\tilde{\chi}_{j(1),...,j(n)}^{i(1),...,i(m)}$ is the sum of the
signs of all edge $3$-colorings of the graph which have the free
edges on the top colored by $i(1),...,i(m)$ (in this order) and
the free edges on the bottom colored by $j(1),...,j(n)$ in this
order.
\end{theorem}


Since the functor $\tilde{F}$ satisfies the Penrose identity and, since, for a planar graph $g\in\hom(0,0)$, $\tilde{F}(g)(1)$ is equal to the number of edge $3$-colorings, we have that this functor generalizes the Penrose invariant \cite{8}.


\section{Binary tree and Eliahou-Kryuchkov conjecture}

From now on, we will restrict ourselves to the study of the
subcategory $\mathbf{PCG}$ where $F=\tilde{F}$.


Recall that a graph is Hamiltonian if there close path the passes
for all the vertices of the graph. There is a well-known theorem
on graph coloring theory due to Whitney \cite{9} that states the following.

\begin{theorem}
If every Hamiltonian planar graph is $4$-colorable then the Four
Color Theorem is true.
\end{theorem}

Let call a morphism $g:1\longrightarrow n$ in $\mathbf{PCG}$
generated only by the generators $\lambda$ and $I$ a descendant
binary $n$-tree, and call a morphism $g:n\longrightarrow 1$ in
$\mathbf{PCG}$ generated only by the generators $y$ and $I$ an
ascendant binary $n$-tree.

When we look at this result in its dual form we have that the Four
Color Theorem is equivalent to the following.

\begin{theorem}
If a morphism $g:1\longrightarrow 1$ is a composition of a
descendant binary $n$-tree with an ascendant binary $n$-tree then
$F(g)$ is non-null.
\end{theorem}

In $\mathbf{CG}$ (or $\mathbf{PCG}$) there is a natural
involution, called the adjoint,
$*:\mathbf{CG}^{op}\rightarrow\mathbf{CG}$ defined by
$\lambda^*=y$, $\cap^*=\cup$ and $x^*=x$ (by definition of
involution we have $y^*=\lambda$, $\cup^*=\cap$, $(f\circ
g)^*=g^*\circ f^*$ and $(f\otimes g)^*=f^*\otimes g^*$).
Geometrically this involution takes the form of a reflection of
the graph in an horizontal line (see the next figure).

$$g = \psdiag{8}{16}{morphism1}
 \longrightarrow g^* =
\psdiag{8}{16}{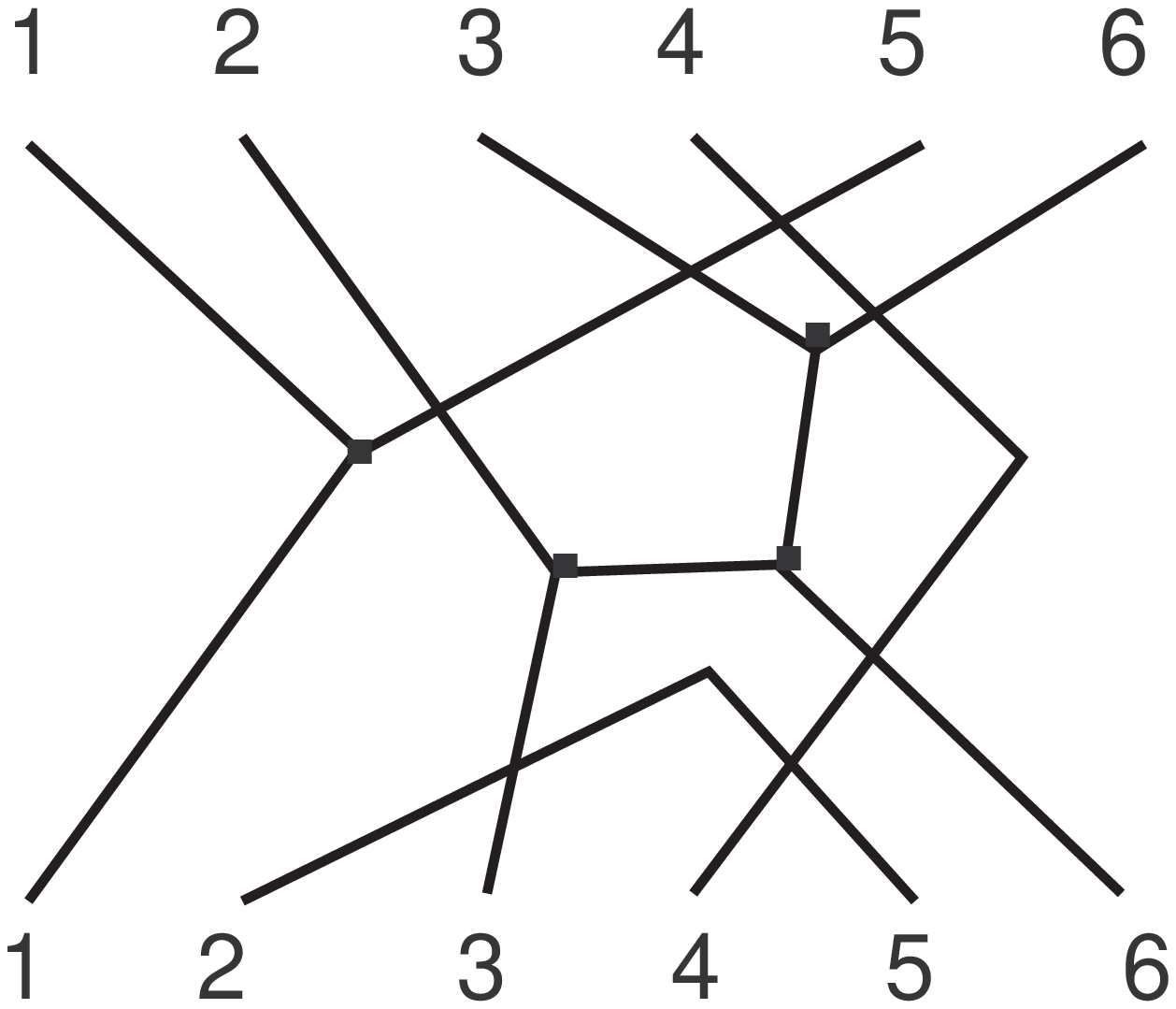}$$

For instance a descendant tree is the adjoint of an ascendant
tree.

Considering the inner products in $\{V^{\otimes
n}\}_{n\in\mathbb{N}}$ defined by their canonical bases and the
involution on $\mathbf{Vect}_\mathbb{K}$ defined by the inner
products ($\langle Tx,y\rangle=\langle x,T^* y\rangle$), we have
that the functors $F$ and $\tilde{F}$ preserve the involution
structure (i.e. $F(g^*)=(F(g))^*$). In particular, we have the
following proposition

\begin{proposition}
If $f:1\rightarrow n$ and $g:1\rightarrow n$ are two descendant
binary $n$-trees then $F(f^*g)=\langle F(g)(e_1),F(f)(e_1)\rangle
id_V$
\end{proposition}

\TeXButton{Proof}{\proof}

We have that $$\begin{array}{cccc}
F(f^*g)(e_i)&=&\sum_{j=1}^{3}\langle F(f^*g)(e_i),e_j \rangle
e_j\\ &=&\sum_{j=1}^{3}\langle F(g)(e_i),F(f)(e_j) \rangle e_j
\end{array}$$

So we have to prove that $\langle F(g)(e_i),F(f)(e_j)
\rangle=\delta_{i,j}\langle F(g)(e_1),F(f)(e_1) \rangle$

If we identify $e_1$, $e_2$ and $e_3$ with the three non-zero
elements of the field $\mathbb{F}_4$ we have that, for any
morphism $g:1\rightarrow n$, if $\langle F(g)(e_i),e_{i(1)}\otimes
\cdots \otimes e_{i(n)} \rangle\not=0$ then
$e_i=e_{i(1)}+e_{i(2)}+ \cdots +e_{i(n)}$.

This proves that $\langle F(g)(e_i),F(f)(e_j) \rangle=0$ if
$i\not=j$.

On the other hand, if $$F(g)(e_i)=\sum
\chi_{i(1),...,i(n)}^{i}e_{i(1)}\otimes \cdots \otimes e_{i(n)}$$
and $\sigma$ is a permutation on $\{1,2,3\}$ then
$$F(g)(\sigma(e_i))=\sum
\chi_{i(1),...,i(n)}^{i}e_{\sigma(i(1))}\otimes \cdots \otimes
e_{\sigma(i(n))}$$

This proves that $\langle F(g)(e_1),F(f)(e_1) \rangle=\langle
F(g)(e_2),F(f)(e_2) \rangle=\langle F(g)(e_3),F(f)(e_3) \rangle$.

\TeXButton{End Proof}{\endproof}

Now we consider the following decomposition of the operator
$F(\lambda)$:

\begin{equation}
F(\lambda)=F(\lambda^+)+F(\lambda^-) \label{eq:+-}
\end{equation}
where $$F(\lambda^+)(e_i)=e_j\otimes e_k \quad \mbox{such
that}\quad (i,j,k)\in\{(1,2,3),(2,3,1),(3,1,2)\}$$ and
$$F(\lambda^-)(e_i)=e_j\otimes e_k \quad \mbox{such that}\quad
(i,j,k)\in\{(2,1,3),(1,3,2),(3,2,1)\}$$

In the same way as a descendant tree is a morphism generated by
$\lambda$ and $I$, a descendant signed tree is a morphism
generated by $\lambda^+$, $\lambda^-$ and $I$.
One simple observation that we can make from (\ref{eq:+-}) is that
for any descendant binary tree $g:1\rightarrow n+1$ with $n$
nodes, $F(g)$ is equal to the sum of the $2^n$ signed trees
corresponding to $g$. Thus $$F(g)(e_1)=\sum
\chi_{i(1),\cdots,i(n+1)}^{1} e_{i(1)}\otimes\cdots \otimes
e_{i(n+1)}$$ with $\chi_{i(1),\cdots,i(n+1)}^{1}=1$ for some $2^n$
indices and zero for the others.

As a consequence of this we have:

\begin{proposition}
For any descendant binary $n$-tree $g:1\rightarrow n+1$ we have
$F(g^*g)=2^n id_V$.
\end{proposition}

Another observation that we can make is a signed reassociation
identity:

$$F((I\otimes\lambda^+)\lambda^+)=F((\lambda^-\otimes I)\lambda^-)
\quad \mbox{and} \quad
F((I\otimes\lambda^-)\lambda^-)=F((\lambda^+\otimes I)\lambda^+)$$

It is well known that any pair of binary trees with the same
number of nodes is connected by a finite sequence of (non-signed)
reassociation moves:

$$\psdiag{10}{20}{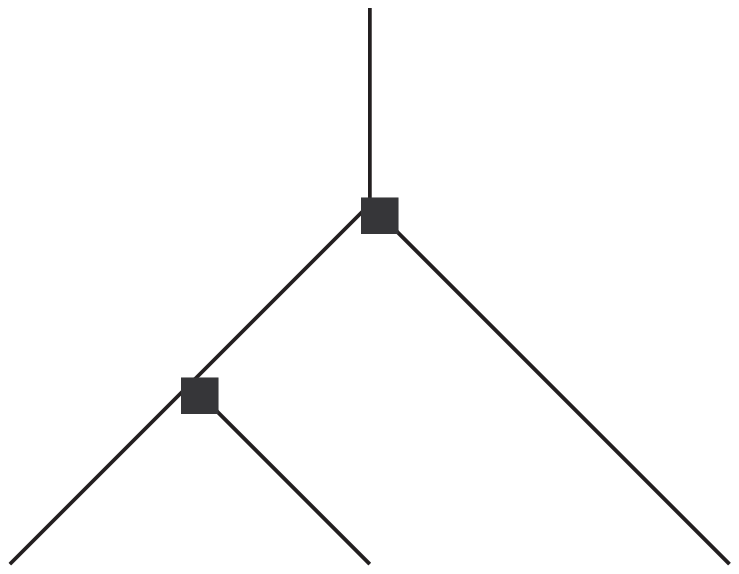}\longleftrightarrow\psdiag{10}{20}{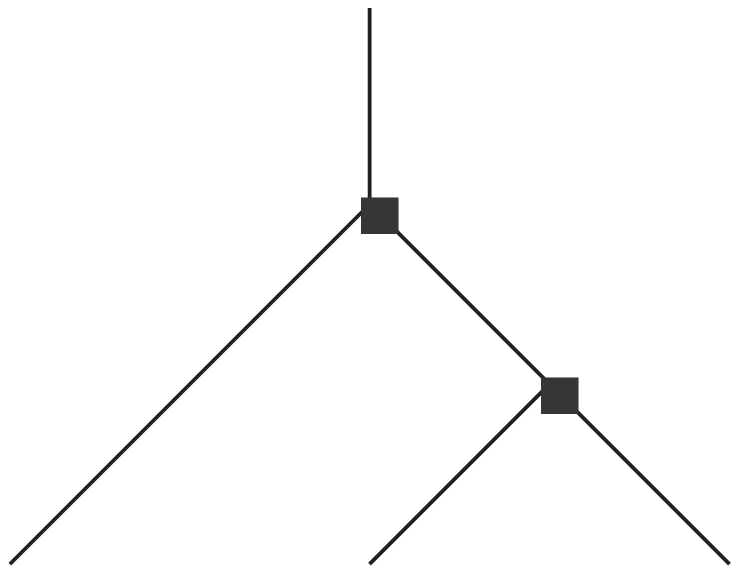}$$

For instance:

$$\psdiag{10}{20}{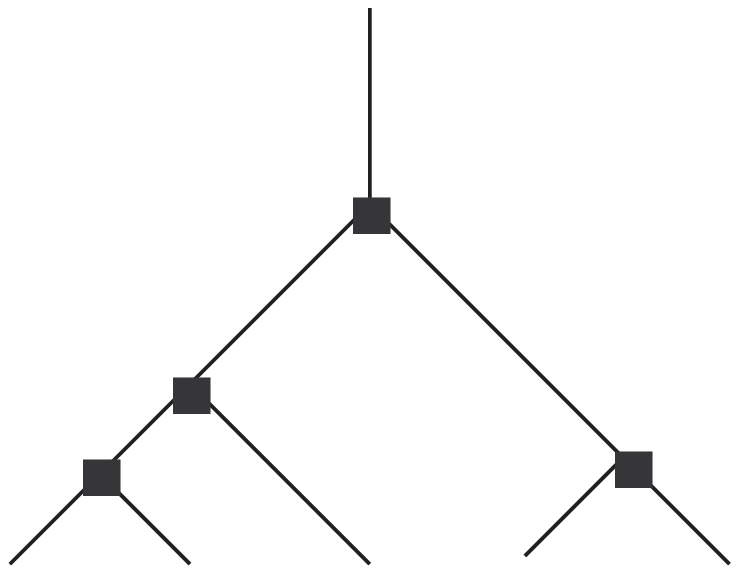}\leftrightarrow\psdiag{10}{20}{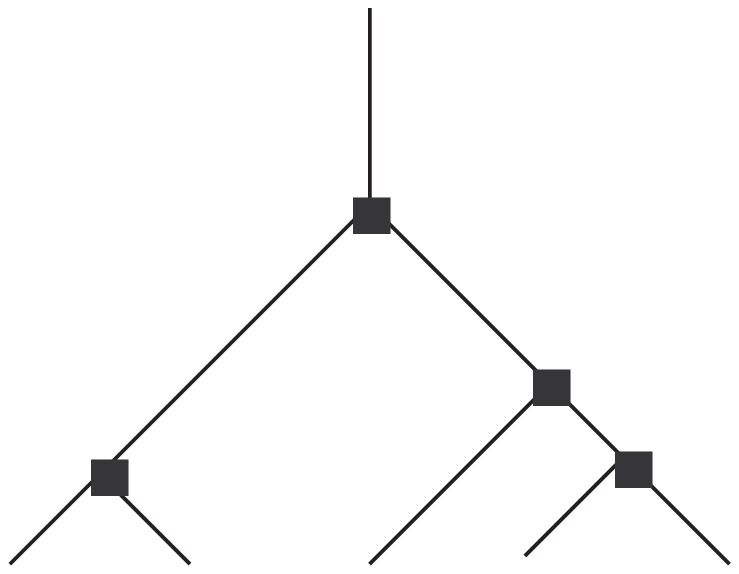}\leftrightarrow
\psdiag{10}{20}{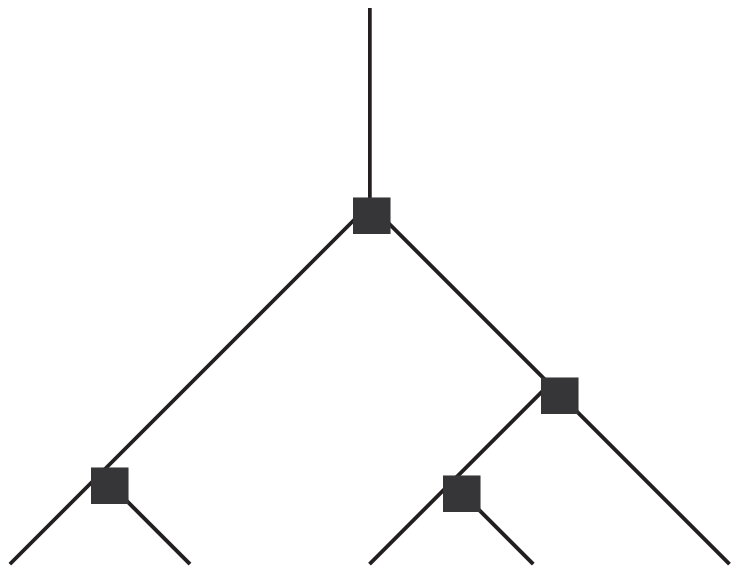}$$

In fact, the graph for which the vertices are the $n$- dimensional
descendant binary trees and the edges represent reassociation
moves between two trees is the $1$-skeleton of the
$n-2$-dimensional {\it associahedron} $\mathcal{A}_{n-2}$ or
Stasheff polytope (see \cite{3} for the definition).

However when we take the analogous graph $\mathcal{A}_{n-2}^{s}$
for signed trees and signed reassociation moves we get a
non-connected graph. Indeed, two signed trees $f$ and $g$ are
connected if $F(f)=F(g)$.

This last graph can be projected in a natural way onto the first
but it is not true that any path on the associahedron
$\mathcal{A}_{n-2}$ can be lifted to a path on
$\mathcal{A}_{n-2}^{s}$.

Eliahou \cite{2} and Kryuchkov (cited from \cite{5}) conjectured the following:

\begin{conjecture}
(Eliahou-Kryuchkov) For any pair of vertices on
$\mathcal{A}_{n-2}$ there exists a path connecting them that can
be lifted to a path on the graph $\mathcal{A}_{n-2}^s$
\end{conjecture}

It is easy to see that this conjecture implies the Four Color
Theorem since two signed trees connected by a sequence of signed
reassociation moves give the same colors on the ends.

In the paper \cite{4} Gravier and Payan proved that this
conjecture is, in fact, equivalent to the Four Color Theorem.

{\bf Acknowledgment} - I wish to thank Roger Picken for his useful suggestions and comments. Supported by Funda\c c\~ao para a Ci\^encia e a Tecnologia, project New Geometry and Topology, PTDC/MAT/101503/2008.

\end{document}